\documentclass[11pt,a4paper]{article}
\usepackage{amsmath}
\usepackage{amsthm}
\usepackage{amssymb}
\usepackage{dsfont}
\usepackage[T1]{fontenc}
\usepackage[french,english,russian]{babel}
\usepackage[utf8]{inputenc}
\usepackage{mathrsfs}
\usepackage{lmodern}
\usepackage{graphicx}
\newtheorem{lemma}{Lemma}
\newtheorem{sublemma}{Sublemma}
\newtheorem{prop}{Proposition}
\newtheorem{theorem}{Theorem}
\newtheorem*{theorems}{Theorem}

\newtheorem{defi}{Definition}
\makeatletter
\def\moverlay{\mathpalette\mov@rlay}
\def\mov@rlay#1#2{\leavevmode\vtop{%
   \baselineskip\z@skip \lineskiplimit-\maxdimen
   \ialign{\hfil$\m@th#1##$\hfil\cr#2\crcr}}}
\newcommand{\charfusion}[3][\mathord]{
    #1{\ifx#1\mathop\vphantom{#2}\fi
        \mathpalette\mov@rlay{#2\cr#3}
      }
    \ifx#1\mathop\expandafter\displaylimits\fi}
\makeatother

\newcommand{\cupdot}{\charfusion[\mathbin]{\cup}{\cdot}}

\title{Limit laws in the lattice problem. \\
V. The case of analytic and stricly convex sets}
\author{Julien Trevisan}
\begin{document}
\selectlanguage{english}
\maketitle
\bigskip
\selectlanguage{french}
\begin{abstract}
Nous étudions l'erreur du nombre de points d'un réseau unimodulaire qui tombent dans un ensemble strictement convexe et analytique possédant l'origine et qui est dilaté d'un facteur $t$. Le but est de généraliser le résultat de $\cite{trevisan2021limit}$. On montre d'abord que l'étude de l'erreur, lorsqu'elle est normalisée par $\sqrt{t}$, lorsque ce paramètre tend vers l'infini et lorsque le réseau considéré est aléatoire, se ramène à l'étude d'une transformée de Siegel $\mathcal{S}(f_{t})(L)$ qui dépend de $t$. Ensuite, on se ramène à l'étude du comportement asymptotic d'une transformée de Siegel avec poids aléatoires, $\mathcal{S}(F)(\theta,L)$ où $\theta$ est un second paramètre aléatoire. Puis, on montre que cette dernière quantité converge presque sûrement et on étudie l'existence des moments de sa loi. Enfin, on montre que ce résultat est encore valable si l'on translate, après dilatation, l'ensemble strictement convexe d'un vecteur $\alpha \in \mathbb{R}^{2}$ fixé. 
\end{abstract}
\selectlanguage{english}
\begin{abstract}
We study the error of the number of points of a unimodular lattice that fall in a strictly convex and analytic set having the origin and that is dilated by a factor $t$. The aim is to generalize the result of $\cite{trevisan2021limit}$. We first show that the study of the error, when it is normalized by $\sqrt{t}$, when this parameter tends to infinity and when the considered lattice is random, is reduced to the study of a Siegel transform $\mathcal{S}(f_{t})(L)$ which depends on $t$. Then, we come back to the study of the asymptotic behaviour of a Siegel transform with random weights, $\mathcal{S}(F)(\theta,L)$ where $\theta$ is a second random parameter. Then, we show that this last quantity converges almost surely and we study the existence of moments of its law. Finally, we show that this result is still valid if we translate, after dilation, the strictly convex set of a fixed vector $\alpha \in \mathbb{R}^{2}$. 
\end{abstract}
\section{Introduction}
The lattice point problem is an open problem in the Geometry of Numbers, at least since Carl Friedrich Gauss took interest in which became the famous Gauss circle problem. The general problem states as followed. \\
Let $d$ be an integer greater than $1$. We recall the following definition : 
\begin{defi}
A subset $L$ of $\mathbb{R}^{d}$ is a lattice if it is a subgroup of $\mathbb{R}^{d}$ such that $L$ is discrete and $\text{span}(L) = \mathbb{R}^{d}$. 
\end{defi}
 Let $P$ be a measurable subset of $\mathbb{R}^{d}$ of non-zero finite Lebesgue measure. We want to evaluate the following cardinal number when $t \rightarrow \infty$ : $$ N(tP + X, L) = | (tP + X) \cap L|$$ where $X \in \mathbb{R}^{d}$, $L$ is a lattice of $\mathbb{R}^{d}$ and $t P + X$ denotes the set $P$ dilated by a factor $t$ relatively to $0$ and then translated by the vector $X$.  \\
Under mild regularity conditions on the set $P$, one can show that : 
$$ N(tP + X, L) = t^{d}\frac{\text{Vol}(P)}{\text{Covol}(L)} + o(t^{d}) $$
where $o(f(t))$ denotes a quantity such that, when divided by $f(t)$, it goes to $0$ when $t \rightarrow \infty$ and where $\text{Covol}(L)$ is defined in the following definition : 
\begin{defi}
The covolume of a lattice $L$ of $\mathbb{R}^{d}$, $\text{covol}(L)$, is the Lebesgue measure of a measurable fundamental set of $L$. Furthermore, a lattice is said to be unimodular if its covolume is equal to $1$. \\
When $d=2$, instead of using the term covolume, we use the term coarea.
\end{defi}
We are interested in the error term $$\mathcal{R}(tP + X,L) = N(tP + X, L) - t^{d}\frac{\text{Vol}(P)}{\text{Covol}(L)} \textit{.}$$
In the case where $d=2$ and where $P$ is the unit disk $\mathbb{D}^{2}$, Hardy's conjecture in $\cite{hardy1917average}$ stipulates that we should have for all $\epsilon > 0$, $$\mathcal{R}(t \mathbb{D}^{2}, \mathbb{Z}^{2}) = O(t^{\frac{1}{2}+\epsilon}) \textit{.}$$
One of the result in this direction has been established by Iwaniec and Mozzochi in $\cite{iwaniec1988divisor}$. They have proven that for all $\epsilon > 0$, $$\mathcal{R}(t \mathbb{D}^{2}, \mathbb{Z}^{2}) = O(t^{\frac{7}{11}+\epsilon}) \textit{.}$$
This result has been recently improved by Huxley in \cite{huxley2003exponential}. Indeed, he has proven that : $$\mathcal{R}(t \mathbb{D}^{2}, \mathbb{Z}^{2}) = O(t^{K} \log(t)^{\Lambda}) $$
where $K = \frac{131}{208} $ and $\Lambda = \frac{18627}{8320}$.  \\
In dimension 3, Heath-Brown has proven in $\cite{heath2012lattice}$ that : 
$$\mathcal{R}(t \mathbb{D}^{3}, \mathbb{Z}^{3}) = O(t^{\frac{21}{16}+\epsilon}) \textit{.}$$
These last two results are all based on estimating what are called $\textit{exponential sums}$.\\ Another approach was followed first by Heath-Brown in $\cite{heath1992distribution}$ and then by Bleher, Cheng, Dyson and Lebowitz in $\cite{bleher1993distribution}$. They took interest in the case where the dilatation parameter $t$ is random. More precisely, they assumed that $t$ was being distributed according to the measure $\rho(\frac{t}{T}) dt$ (that is absolutely continuous relatively to Lebesgue measure) and where $\rho$ is a probability density on $[0,1]$ and $T$ is parameter that goes to infinity. In that case, Bleher, Cheng, Dyson and Lebowitz showed the following result (which generalizes the result of Heath-Brown) : 
\begin{theorems}[$\cite{bleher1993distribution}$]
Let $\alpha \in [0,1[^{2}$. There exists a probability density $p_{\alpha}$ on $\mathbb{R}$ such that for every piecewise continuous and bounded function $g : \mathbb{R} \longrightarrow \mathbb{R}$, 
$$ \lim_{T \rightarrow \infty} \frac{1}{T} \int_{0}^{T} g(\frac{\mathcal{R}(t \mathbb{D}^{2} + \alpha , \mathbb{Z}^{2})}{\sqrt{t}}) \rho(\frac{t}{T}) dt = \int_{\mathbb{R}} g(x) p_{\alpha}(x) dx \textit{.}$$
Furthermore $p_{\alpha}$ can be extended as an analytic function over $\mathbb{C}$ and verifies that for every $\epsilon > 0$, $$p_{\alpha}(x) = O(e^{-|x|^{4- \epsilon}})$$ when $x \in \mathbb{R}$ and when $|x| \rightarrow \infty$. 
\end{theorems}
In our case, we keep $t$ deterministic as in the original Gauss problem but we let the lattice $L$ be a random unimodular lattice and we study $\mathcal{R}$. This approach was first initiated by Kesten in $\cite{kesten60u}$ and in $\cite{kesten62u}$. It should be noted that several counting problems have followed this approach : we can cite, for example, $\cite{gorodnikcounting}$, $\cite{dolgopyat2014deviations}$, $\cite{marklof1998n}$, $\cite{beck1994probabilistic}$, $\cite{levin2013gaussian}$ and $\cite{vinogradov2010limiting}$. \\
We denote by $\mathscr{S}_{2}$ the space of unimodular lattices and it can be seen as the quotient space $SL_{2}(\mathbb{R})/SL_{2}(\mathbb{Z})$. We denote by $\mu_{2}$ the unique Haar probability measure on it. Let us set :
\begin{equation}
\label{eq67}
\Pi = \{ (k_{1},k_{2}) \in \mathbb{Z}^{2} \text{ } | \text{ } k_{1} \wedge k_{2} = 1 \textit{, } k_{1} \geqslant 0  \}
\end{equation}
where we agree that if $k_{1} = 0$, $k_{1} \wedge k_{2} = 1$ means that $k_{2} = 1$. \\
We denote by $\lVert \cdot \rVert$ the usual euclidean norm over $\mathbb{R}^{2}$.\\
We need to define some additional objects.
\begin{defi}
\label{def3}
\label{def1}
For every $ i \in \{ 1, 2 \}$, we call $$\lVert L \rVert _{i} = \min \{ r > 0 \text{ } | \text{ } B_{f}(0,r) \text{ contains } i \text{ vectors of L lineary independent} \}  $$
where $B_{f}(0,r)$ is the closed centred ball on $0$ for the norm $\lVert \cdot \rVert$ of radius $r$. \\
\emph{In fact, for almost all $L \in \mathscr{S}_{2}$, $\lVert L \rVert_{2} > \lVert L \rVert_{1} $ and there exists only one couple of vectors $(e_{1}(L),e_{2}(L))$ such that $(e_{1}(L))_{1} > 0$, $ \lVert e_{1}(L) \rVert = \lVert L \rVert_{1}$, $(e_{2}(L))_{1} > 0$ and $\lVert e_{2}(L) \rVert = \lVert L \rVert_{2}$.} \\
\emph{In the rest of the article, for the sake of simplicity, for $L$ a lattice, we will use the notation $\lVert L \rVert$ instead of $\lVert L \rVert_{1}$.} \\
For a lattice $L \in \mathscr{S}_{2}$, we also say that a vector of $L$ is prime if it is not a non-trivial integer multiple of another vector of $L$. \\
\emph{In fact, for every $M \in SL_{2}(\mathbb{R})$,  a vector $ l \in L$ is prime if, and only if, $M l \in M L$ is prime and a vector $(k_{1},k_{2}) \in \mathbb{Z}^{2}$ is prime if, and only if, $k_{1} \wedge k_{2} = 1$. With these notations, one has that, for a generic lattice $L \in \mathscr{S}_{2}$, $ e \in L$ is prime if, and only if, $e$ can be written as $ e=k_{1} e_{1}(L) + k_{2} e_{2}(L)$ with $k_{1} \wedge k_{2} = 1$.} \\
Finally, for a generic lattice $L \in \mathscr{S}_{2}$, we call $\mathcal{P}_{+}(L)$ the set of vectors $e$ of $L$ such that $e=k_{1} e_{1}(L) + k_{2} e_{2}(L)$ with $(k_{1},k_{2}) \in \Pi$. \emph{All the vectors of $\mathcal{P}_{+}(L)$ are prime vectors according to the previous remark.}
\end{defi}
We recall also the fact that we say that a real random variable $Z$ is symmetrical if $\mathbb{P}_{Z} = \mathbb{P}_{-Z} $ where, for every random variable $X$, $\mathbb{P}_{X}$ stands for the law of the random variable $X$.\\
Let $\tilde{\mu}_{2}$ be a probability measure that has a smooth bounded density $\sigma$ with respect to $\mu_{2}$. There are two different cases that are addressed in our main result, which is Theorem $\ref{thm2}$. The first one is when $\tilde{\mu}_{2}$ is compactly supported, $\textit{id est}$ when there exists $\alpha > 0$ such that : 
\begin{equation}
\label{eq1000}
\tilde{\mu}_{2}(\{ L \in \mathscr{S}_{2} \text{ } | \text{ } \lVert L \rVert_{1} < \alpha \}) = 0 \textit{.}
\end{equation}
The second one is when $\tilde{\mu}_{2}$ is non-compactly supported under the following condition : there exists $m > 0$, there exists $\alpha > 0$ such that for all $L$ that belongs to the event $( \lVert L \rVert < \alpha )$,
\begin{equation}
\label{eq1001}
\sigma(L) \geqslant m \textit{.}
\end{equation}
An example of such a measure $\tilde{\mu}_{2}$ is given by the normalized Haar measure $\mu_{2}$. \\
Let $\mathcal{E}$ be an ellipse centred around $0$. Let's call $M$ a matrix that transforms $\mathcal{E}$ into a disk and that belongs to $SL_{2}(\mathbb{R})$. $M$ is unique modulo the natural action of $SL_{2}(\mathbb{Z})$. \\
Let's set : $\mathbb{T}^{\infty} = (\mathbb{T}^{1})^{\Pi}$
where $\mathbb{T}^{1} = \mathbb{R}/\mathbb{Z}$ and let's call $\lambda_{\infty}$ the normalized Lebesgue measure product over $\mathbb{T}^{\infty}$. \\
The main result of the previous article $\cite{trevisan2021limit}$ is the following theorem :  
\begin{theorem}
\label{thm2}
For every real numbers $a <b$, 
$$\lim_{t \rightarrow \infty} \tilde{\mu}_{2} \left( L \in \mathscr{S}_{2} \textit{ }| \textit{ } \frac{\mathcal{R}(t \mathcal{E}, L)}{\sqrt{t}} \in [a,b] \right) = (\lambda_{\infty} \times (M^{-1})_{*}\tilde{\mu}_{2} ) \left( (\theta,L) \in \mathbb{T}^{\infty} \times \mathscr{S}_{2} \textit{ }| \textit{ }  S(\theta,L^{\perp}) \in [a,b] \right) $$
where $(M^{-1})_{*}\tilde{\mu}_{2}$ is the push-forward of $\tilde{\mu}_{2}$ by $M^{-1}$ and where $\theta = (\theta_{e}) \in \mathbb{T}^{\infty}$,
$$S(\theta,L) = \frac{2}{\pi} \sum_{e \in  \mathcal{P}_{+}(L) } \frac{\phi(\theta_{e})}{\lVert e \rVert^{\frac{3}{2}}}  $$ 
with 
\begin{equation}
\label{eq501}
\phi(\theta_{e})= \sum_{m \geqslant 1} \frac{\cos(2 \pi m \theta_{e} - \frac{3 \pi}{4})}{m^{\frac{3}{2}}} \textit{.}
\end{equation}
Furthermore, $S(\theta,L)$ (and $S(\theta,L^{\perp}$) : \begin{itemize}
\item converges almost surely
\item is symmetrical and its expectation is equal to $0$
\item admits a moment of order $1+\kappa$ for any $0 \leqslant \kappa < \frac{1}{3}$ 
\item $S(\theta,L)$ admits moments of all order $1 \leqslant p < \infty$ when $\tilde{\mu}_{2}$ is compactly supported
\item does not admit a moment of order $\frac{4}{3}$ when $\tilde{\mu}_{2}$ is non-compactly supported under the condition ($\ref{eq1001}$). 
\end{itemize}
\end{theorem}
In this article, we want to extend this last result. More precisely, let's suppose that $\mathcal{E} = \Omega_{\gamma}$ where $\gamma$ is analytic curve that is simple, closed and strictly convex and where $0 \in \Omega_{\gamma}$. Let's also call $x_{\gamma}(\xi)$ the point on $\gamma$ where the outer normal to $\gamma$ coincides with $\frac{\gamma}{\lVert \gamma \rVert}$ and $\rho_{\gamma}(\xi)$ the curvature radius of $\gamma$ at $x_{\gamma}(\xi)$. Let's set finally $\mathbb{T}^{\infty,2} = (\mathbb{T}^{1})^{\Pi} \times (\mathbb{T}^{1})^{\Pi} $
and let's call $\lambda_{\infty,2}$ the normalized Lebesgue measure product over $\mathbb{T}^{\infty,2}$. \\ 
Then, we want to prove the following theorem :
\begin{theorem}
\label{thm3}
There exists a distribution function $\mathcal{D}_{\gamma}(z)$ such that for every real $z$ we have : 
$$\lim_{t \rightarrow \infty} \tilde{\mu}_{2} \left( \frac{\mathcal{R}(t \Omega_{\gamma}, L)}{\sqrt{t}} \in ]- \infty, z] \right) = \mathcal{D}_{\gamma}(z) \textit{.} $$
In the case where $\Omega_{\gamma}$ is symmetric, $\mathcal{D}_{\gamma}(z)$ is the distribution function of 
$$ S_{\gamma}(\theta,L) = \frac{2}{\pi} \sum_{e \in  \mathcal{P}_{+}(L) } \frac{\rho_{\gamma}(e) \phi(\theta_{e})}{\lVert e \rVert^{\frac{3}{2}}} $$
with $\theta = (\theta_{e}) \in \mathbb{T}^{\infty}$ being distributed according to $\lambda_{\infty}$ and $L$ being distributed according to $\tilde{\mu}_{2}$. \\
In the non symmetric case,  $\mathcal{D}_{\gamma}(z)$ is the distribution function of 
$$ S_{\gamma}(\theta,L) = \frac{1}{\pi} \sum_{e \in  \mathcal{P}_{+}(L) } \frac{ \phi_{\gamma,2}(\theta_{e},e)}{\lVert e \rVert^{\frac{3}{2}}} $$
where $\theta = (\theta_{1,e}, \theta_{2,e}) \in \mathbb{T}^{\infty,2}$ being distributed according to $\lambda_{\infty,2}$ and $L$ being distributed according to $\tilde{\mu}_{2}$ and with
\begin{equation}
\label{eq0501}
\phi_{\gamma,2}(\theta_{e},e)= \sum_{m \geqslant 1} \frac{\rho_{\gamma}(e)\cos(2 \pi m \theta_{1,e} - \frac{3 \pi}{4}) + \rho_{\gamma}(-e)\cos(2 \pi m \theta_{2,e} - \frac{3 \pi}{4}) }{m^{\frac{3}{2}}} \textit{.}
\end{equation}
Furthermore, in both cases, $S_{\gamma}$ admits the same properties of $S(\theta,L)$ listed in Theorem $\ref{thm2}$.
\end{theorem}
In fact, this theorem can be generalized as followed : 
\begin{theorem}
\label{thm4}
For every $\alpha \in \mathbb{R}^{2}$, there exists a distribution function $\mathcal{D}_{\gamma,\alpha}(z)$ such that for every real $z$ we have : 
$$\lim_{t \rightarrow \infty} \tilde{\mu}_{2} \left( \frac{\mathcal{R}(t \Omega_{\gamma}+\alpha, L)}{\sqrt{t}} \in ]- \infty, z] \right) = \mathcal{D}_{\gamma,\alpha}(z) \textit{.} $$
In the case where $\Omega_{\gamma}$ is symmetric, $\mathcal{D}_{\gamma}(z)$ is the distribution function of 
$$ S_{\gamma}(\theta,L) = \frac{2}{\pi} \sum_{e \in  \mathcal{P}_{+}(L) } \frac{\rho_{\gamma}(e) \phi_{\alpha}(\theta_{e},e)}{\lVert e \rVert^{\frac{3}{2}}} $$
with $\theta = (\theta_{e}) \in \mathbb{T}^{\infty}$ being distributed according to $\lambda_{\infty}$ and $L$ being distributed according to $\tilde{\mu}_{2}$ and where 
$$\phi_{\alpha}(\theta_{e},e)= \sum_{m \geqslant 1} \frac{\cos(2 \pi m \theta_{e} + 2 \pi m <\alpha,e> - \frac{3 \pi}{4})}{m^{\frac{3}{2}}} \textit{.}$$
In the non symmetric case,  $\mathcal{D}_{\gamma,\alpha}(z)$ is the distribution function of 
$$ S_{\gamma}(\theta,L) = \frac{1}{\pi} \sum_{e \in  \mathcal{P}_{+}(L) } \frac{ \phi_{\alpha,\gamma,2}(\theta_{e},e)}{\lVert e \rVert^{\frac{3}{2}}} $$
where $\theta = (\theta_{1,e}, \theta_{2,e}) \in \mathbb{T}^{\infty,2}$ being distributed according to $\lambda_{\infty,2}$ and $L$ being distributed according to $\tilde{\mu}_{2}$ and with
$$\phi_{\alpha, \gamma,2}(\theta_{e},e)= \sum_{m \geqslant 1} \frac{\rho_{\gamma}(e)\cos(2 \pi m \theta_{1,e} + 2 \pi m <\alpha,e> - \frac{3 \pi}{4}) + \rho_{\gamma}(-e)\cos(2 \pi m \theta_{2,e} - 2 \pi m <\alpha,e> - \frac{3 \pi}{4}) }{m^{\frac{3}{2}}} \textit{.}$$
Furthermore, in both cases, $S_{\gamma}$ admits the same properties of $S(\theta,L)$ listed in Theorem $\ref{thm2}$.
\end{theorem}
This last theorem is a generalization of Theorem $\ref{thm2}$ on two planes. The first one is the shape of the sets : it treats the more general case of analytic curves, not only the case of the ellipses centred on $0$. The second one is the presence of a translation parameter $\alpha$ whereas Theorem $\ref{thm2}$ can be deduced by assuming that $\alpha = 0$ (no translation of $t \mathcal{E}$). \\
In the next section, we give a brief heuristic explanation of the approach that we will follow. It is basically the same as in $\cite{trevisan2021limit}$. At the end of the section, we will give the plan of the paper.
\section{Heuristic explanation and plan of the proof}
First, let's explain the different steps of the proof of Theorem $\ref{thm3}$, Theorem $\ref{thm4}$ being a simple generalization of Theorem $\ref{thm3}$. \\
$\textit{First step.}$ By regularizing the problem and using the Poisson summation formula, we are going to show that the quantity $S_{A,prime}(L,t)$, for $A > 0$ a fixed parameter that is taken large enough, when $t \rightarrow \infty$, is close, in probability, to $\frac{\mathcal{R}(t \Omega_{\gamma},L)}{\sqrt{t}}$ where $S_{A,prime}(L,t)$ is defined by : 
\begin{equation}
\label{eq44}
S_{A,prime}(L,t) = \frac{1}{\pi} \sum_{\substack{ l \in L^{\perp} \textit{ prime}  \\ 0 < \lVert l \rVert \leqslant A }} \frac{ 1}{\lVert l \rVert^{\frac{3}{2}}} \sum_{m \in \mathbb{N}-\{ 0 \}} \frac{ \rho_{\gamma}(l) \cos(2 \pi t m Y_{\gamma}(l) - \frac{3 \pi}{4}) + \rho_{\gamma}(-l) \cos(2 \pi t m Y_{\gamma}(-l) - \frac{3 \pi}{4})  }{m^{\frac{3}{2}}}
\end{equation}
where $Y_{\gamma}(\xi)$, for all $\xi \in \mathbb{R}^{2}-\{ 0 \}$, is defined by 
$$Y_{\gamma}(\xi) = < \xi, x_{\gamma}(\xi) > $$
(it is a positive homogeneous function of order 1). \\
Let's say a few remarks. First, the vectors $l$ prime considered are such that $\lVert l \rVert \leqslant A$. We have to limit the norm of the considered vectors because of a convergence problem. Second, we take into account a phenomenon of multiplicity (if $l$ appears in the sum, $2l$ is also going to appear). In $\cite{bleher1993distribution}$ such a phenomenon was also taken into account. This was done in order to get independence at infinity, as it was also done in $\cite{bassam}$, and we do that for the same goal.  \\
By using the remark that is in Definition $\ref{def3}$ and by replacing $L$ by $L^{\perp}$ (which is done only for a matter of convenience), one has that : 
\begin{align}
\label{eq200}
S_{A,prime}(L^{\perp},t) = \frac{1}{\pi} \sum_{\substack {k_{1} \wedge k_{2} = 1 \\ k_{1} \geqslant 0 \\ \lVert k_{1} e_{1}(L) + k_{2} e_{2}(L) \rVert \leqslant A }} \frac{1}{\lVert k_{1} e_{1}(L) + k_{2} e_{2}(L) \rVert^{\frac{3}{2}}}  & \\
\phi_{\gamma,2}\left((tY_{\gamma}(k_{1} e_{1}(L) + k_{2} e_{2}(L)),tY_{\gamma}(-(k_{1} e_{1}(L) + k_{2} e_{2}(L)))),k_{1} e_{1}(L) + k_{2} e_{2}(L) \right) &  \nonumber
\end{align}
where the function $\phi_{\gamma,2}$ was defined by the Equation $(\ref{eq0501})$. We have done that so from this stage onwards we consider vectors of $\mathcal{P}_{+}(L)$ with a fixed indexation (that does not depend on $L$). Furthermore, this indexation will be very useful for the second step (for more details, see Section 4). \\
\\
$\textit{Second step.}$ 
In the non symmetric case, we will show that the family of variables $(tY_{\gamma}(k_{1} e_{1}(L) + k_{2} e_{2}(L)),tY_{\gamma}(-(k_{1} e_{1}(L) + k_{2} e_{2}(L))))$, whose values are in $(\mathbb{R}/ \mathbb{Z})^{2}$, become, when $t \rightarrow \infty$, independent from one another and indeed converge towards independent and identically distributed random variables whose common distribution is given by the normalized Haar measure on $(\mathbb{R}/ \mathbb{Z})^{2}$. The idea here is basically the same as in $\cite{bleher1993distribution}$ and in $\cite{heath1992distribution}$ where the respective authors used the fact that the square roots of square free integers are $\mathbb{Z}$-free. It is a generalization of what was done in $\cite{trevisan2021limit}$. In our case, to prove the result, we will decompose the space of unimodular lattices into small geodesic segments, calculate the Taylor series of $Y_{\gamma}(k_{1} e_{1}(L) + k_{2} e_{2}(L))$ and of $Y_{\gamma}(-(k_{1} e_{1}(L) + k_{2} e_{2}(L)))$ at order 1 on such a segment and show that the coefficients of order 1 are $\mathbb{Z}$-free. \\
We will also prove that these variables become independent, when $t \rightarrow \infty$, from the variable $L$ due to the presence of the factor $t$. \\
In the symmetric case, instead of considering  $(tY_{\gamma}(k_{1} e_{1}(L) + k_{2} e_{2}(L)),tY_{\gamma}(-(k_{1} e_{1}(L) + k_{2} e_{2}(L))))$, we consider the family of variables $((tY_{\gamma}(k_{1} e_{1}(L) + k_{2} e_{2}(L)))$, whose values are in $(\mathbb{R}/ \mathbb{Z})$, and show the same results.
\\
\\
$\textit{Third step.}$
Thanks to the first and second step, we will see that the asymptotic distribution of $\frac{\mathcal{R}(t \Omega_{\gamma},L)}{\sqrt{t}}$ is the distribution of $S_{\gamma}(\theta,L)$ (see Theorem $\ref{thm3}$) under the assumption that the quantity $S(\theta,L)$ is well-defined. This last fact will be quasi immediate because what was done in the last section of $\cite{trevisan2021limit}$ can be generalized directly in our case. Furthermore, all the listed properties of $S_{\gamma}(\theta,L)$ are also going to be obtained immediately. \\
After doing all of that, we will finally get the validity of Theorem $\ref{thm3}$. \\
\\
$\textbf{Plan of the paper.}$
The next section will be dedicated to deal with the first step of the proof, namely it will show that $\frac{\mathcal{R}(t \Omega_{\gamma},L)}{\sqrt{t}}$ is close in probability with $S_{A,prime}(L^{\perp},t) $
when $A$ is a fixed parameter taken large enough and $t$ goes to infinity (see Proposition $\ref{prop8}$). We have to "cut" the sum because of the problem of convergence of the Fourier series of $X  \longmapsto \frac{\mathcal{R}(t \Omega_{\gamma}+X,L)}{\sqrt{t}}$ which is due to the lack of regularity of the indicator function $\mathbf{1}_{t \Omega_{\gamma}}$. To prove this, we are going to proceed by $\textit{regularization}$ which means here that we are going to smooth the indicator function $\mathbf{1}_{t \Omega_{\gamma}}$ via a Gaussian kernel. \\
In Section 4, we tackle the second step of the proof, that is the fact, in the non symmetric case, that the $(t Y_{\gamma}(k_{1}e_{1}(L)+k_{2}e_{2}(L),t Y_{\gamma}(-(k_{1}e_{1}(L)+k_{2}e_{2}(L)))$ become independent when $t \rightarrow \infty$. We also show that they converge towards random variables that are identically distributed according to the normalized Haar measure over $(\mathbb{R}/\mathbb{Z})^{2}$ and that they become independent, when $t \rightarrow \infty$, from $L$ and so that $\frac{\mathcal{R}(t \Omega_{\gamma},L)}{\sqrt{t}}$ has the same distribution of $S_{\gamma}(\theta,L)$. We also deal with the symmetric case, which is simpler. \\
In Section 5 we are going to tackle the third step of the proof, namely study the convergence of $\tilde{S}_{A}(\omega, L)$ when $A \rightarrow \infty$ and the existence of moments of its limit. \\
In Section 6, which is the last section, we give the approach, based on the approach to prove Theorem $\ref{thm3}$, to prove Theorem $\ref{thm4}$. \\
In the rest of the article, all the calculus of expectation $\mathbb{E}$, of variance $\text{Var}$ and of probability $\mathbb{P}$ will be made according to the measure $\tilde{\mu}_{2}$. Furthermore, the expression typical is going to signify $\tilde{\mu}_{2}-\textit{almost surely}$. In fact, like we have said in Section 1, we are going to suppose that $\tilde{\mu}_{2} = \mu_{2}$ in Section 3 and in Section 4 because all the results extend to the general case.
\section{Reduction to the study of the Siegel transform}
The main object of this section is to show the following proposition : 
\begin{prop} 
\label{prop8}
For every $\alpha > 0$, for every $A > 0$ large enough, for every $t$ large enough, one has that :
$$\mathbb{P}(\Delta_{A,prime}(L,t)  \geqslant \alpha) \leqslant \alpha$$
where 
\begin{equation}
\label{eq201}
\Delta_{A,prime}(L,t) = |\frac{\mathcal{R}(t \Omega_{\gamma},L)}{\sqrt{t}} - S_{A,prime}(L,t)| \textit{ .}
\end{equation}

\end{prop}
This proposition basically says that we can reduce the asymptotical study of $\frac{\mathcal{R}(t \Omega_{\gamma},L)}{\sqrt{t}}$ to the study of its Fourier transform, taking into account a phenomenon of multiplicity. \\
In fact, due to the triangle inequality, we only have to prove the following two lemmas : 
\begin{lemma}
\label{lemme30}
For every $\alpha > 0$, for every $A > 0$ large enough, for every $t$ large enough, one has that :
$$\mathbb{P}(\Delta_{A}(L,t) \geqslant \alpha) \leqslant \alpha$$
where 
\begin{equation}
\label{eq202}
\Delta_{A}(L,t) = |\frac{\mathcal{R}(t \Omega_{\gamma},L)}{\sqrt{t}} - H_{A}(L,t)| \textit{ with}
\end{equation}

\begin{equation}
\label{eq10}
H_{A}(L,t) = \frac{1}{2 \pi} \sum_{\substack{ l \in L^{\perp} \\ 0 < \lVert l \rVert \leqslant A }} \frac{1}{\lVert l \rVert^{\frac{3}{2}} } \left( \rho_{\gamma}(l) \cos(2 \pi t Y_{\gamma}(l) - \frac{3 \pi}{4} )
 +   \rho_{\gamma}(-l) \cos(2 \pi t Y_{\gamma}(-l) - \frac{3 \pi}{4}) \right)
\end{equation}
which becomes when $\Omega_{\gamma}$ is symmetric 
\begin{equation}
\label{eq100000}
H_{A}(L,t) = \frac{1}{ \pi} \sum_{\substack{ l \in L^{\perp} \\ 0 < \lVert l \rVert \leqslant A }} \frac{\rho_{\gamma}(l) \cos(2 \pi t Y_{\gamma}(l) - \frac{3 \pi}{4}))}{\lVert l \rVert^{\frac{3}{2}} }   \textit{.}
\end{equation}
\end{lemma}
\begin{lemma}
\label{lemme31}
For every $\alpha > 0$, for every $A > 0$ large enough, for every $t$ large enough, one has that :
$$\mathbb{P} (|S_{A,prime}(L,t) - H_{A}(L,t)| \geqslant \alpha) \leqslant \alpha \textit{.}$$
\end{lemma}
\begin{proof}[Proof of Proposition $\ref{prop8}$]
One has that : 
\begin{equation}
\label{eq54}
\Delta_{A,prime}(L,t) \leqslant \Delta_{A}(L,t) + |S_{A,prime}(L,t) - H_{A}(L,t)| \textit{.}
\end{equation}
The Lemma $\ref{lemme30}$ and the Lemma $\ref{lemme31}$ imply then the wanted result.
\end{proof}
Let's say a few words about Lemmas $\ref{lemme30}$ and $\ref{lemme31}$ before following with their respective proofs. The Lemma $\ref{lemme30}$ says that the study of $\frac{\mathcal{R}(t \Omega_{\gamma},L)}{\sqrt{t}}$ can be reduced to the study of its Fourier transform. The Lemma $\ref{lemme31}$ says that the phenomenon of multiplicity (the fact that for a prime vector $l$, $2 l$, $3 l$ etc. appear in the sum $H_{A}(L,t)$ when $A \rightarrow \infty$) is not so important. We only have to gather all the multiples of a prime vector (which corresponds to the infinite sum over $m$, see equation ($\ref{eq44}$)), so that we focus on prime vectors. 
\subsection{Proof of Lemma $\ref{lemme30}$}
First, we are going to prove the Lemma $\ref{lemme30}$. To do so, we are following closely the approach of $\cite{bleher1992distribution}$, yet with some differences because in our case it is not the radius of dilatation that is random but the lattice (or, equivalently and in a certain sense, the oval). \\
For $x \in \mathbb{R}^{2}$ and $t > 0$, let's define $$\lambda(x ; t) = \frac{t^{2}}{4 \pi} e^{-\frac{t^{2}}{4 \pi} \lVert x \rVert^{2}}$$
and, for $M \in SL_{2}(\mathbb{R})$,
\begin{equation}
\label{eq11}
\lambda_{M}(x;t)=\lambda(M x ; t) \textit{.}
\end{equation}
We recall that : 
\begin{equation}
\label{eq12}
\int_{\mathbb{R}^{2}} \lambda_{M}(x;t) dx = 1 
\end{equation}
and that the Fourier transform of $\lambda_{M}(\cdot; t)$ can be expressed as 
\begin{equation}
\label{eq13}
 \widetilde{\lambda_{M}}(\xi ; t) = e^{- \frac{\lVert (M^{-1})^{T} \xi \rVert^{2}}{t^{2}}} \textit{.}
\end{equation}
We introduce the following function : 
\begin{equation}
\label{eq14}
\chi_{\gamma, M}(x;t) = (\mathbf{1}_{t \Omega_{M^{-1}\gamma}}*\lambda_{M}(\cdot; t))(x) = \int_{\mathbb{R}^{2}} \mathbf{1}_{t \Omega_{M^{-1}\gamma}}(y) \lambda_{M}(x-y; t) dy 
\end{equation}
(it is a regularization of the function $\mathbf{1}_{t \Omega_{M^{-1}\gamma}}$). \\
Let us also set : 
\begin{equation}
\label{eq15}
N_{reg}(t \Omega_{\gamma}, M) = \sum_{n \in \mathbb{Z}^{2}} \chi_{\gamma, M}(n;t) \textit{ and }
\end{equation}
(the index "reg" stands for regularized) 
\begin{equation}
\label{eq16}
F(M,t) = \frac{N_{reg}(t \Omega_{\gamma}, M) - \text{Area}(t \Omega_{\gamma})}{\sqrt{t}} \textit{.}
\end{equation}
Let $L$ be a unimodular lattice such that $e_{1}(L)$ and $e_{2}(L)$ are well-defined and let 
\begin{equation}
\label{eq17}
M=[e_{1}(L),e_{2}(L)] \textit{ if } \text{det}([e_{1}(L),e_{2}(L)]) > 0 
\end{equation}
and
\begin{equation}
\label{eq18}
M=[e_{2}(L),e_{1}(L)] \textit{ if } \text{det}([e_{2}(L),e_{1}(L)]) > 0 \textit{.}
\end{equation}
Then $M$ is a matrix that represents $L$ and one has immediately that  : 
$$\mathcal{R}(t \Omega_{\gamma},L)= \mathcal{R}(t \Omega_{M^{-1}\gamma},\mathbb{Z}^{2}) \textit{.}$$
Now, let's call : 
\begin{equation}
\label{eq19}
\Delta_{1}(L,t) = | \frac{\mathcal{R}(t \Omega_{\gamma},L)}{\sqrt{t}} - F(M,t) | \textit{ and } (\Delta_{2})_{A}(L,t) = | F(M,t) - H_{A}(L,t) |
\end{equation}
so one has that : 
\begin{equation}
\label{eq20}
\Delta_{A}(L,t) \leqslant \Delta_{1}(L,t) + (\Delta_{2})_{A}(L,t)\textit{.}
\end{equation}
The proof of Lemma $\ref{lemme30}$ lies on the two following lemmas : 
\begin{lemma}
\label{lemme32}
The quantity $\Delta_{1}(L,t)$ converges to $0$ when $t \rightarrow \infty$. 
\end{lemma}
\begin{lemma}
\label{lemme33}
For all $\alpha > 0$, for all $A$ large enough, for all $t$ large enough,
$$\mathbb{P}((\Delta_{2})_{A}(L,t) \geqslant \alpha) \leqslant \alpha \textit{.}$$
\end{lemma}
\begin{proof}[Proof of Lemma $\ref{lemme30}$]
It is the direct consequence of Equation ($\ref{eq20}$) and of Lemma $\ref{lemme32}$ and Lemma $\ref{lemme33}$.
\end{proof}
Lemma $\ref{lemme32}$ basically tells us that the study of $\frac{\mathcal{R}(t \Omega_{\gamma},L)}{\sqrt{t}}$ can be reduced to the study of one of its regularized Fourier series, whereas Lemma $\ref{lemme33}$ means that the asymptotical study of this regularized Fourier series can be brought back to the study of the non-regularized Fourier series. \\
The next subsubsection is dedicated to the proof of Lemma $\ref{lemme32}$ and the subsubsection after it is dedicated to the proof of Lemma $\ref{lemme33}$. 
\subsubsection{Proof of Lemma $\ref{lemme32}$}
The proof of Lemma $\ref{lemme32}$ is based on two sublemmas. The first one is the following :
\begin{sublemma}
\label{lemme3}
For all $x \in \mathbb{R}^{2}$, for all $t > 0$, $$|\chi_{\gamma, M}(x;t) - \mathbf{1}_{t \Omega_{M^{-1}\gamma}}(x)| \leqslant e^{- \frac{t^{2}}{4} \text{dist}(M x,t \gamma)^{2}} $$ 
where for all $z \in \mathbb{R}^{2}$, $$\text{dist}(z,t \gamma)= \inf_{y \in t \gamma} |z-y| \textit{.}$$
\end{sublemma}
\begin{proof}
One has that : 
$$|\chi_{\gamma, M}(x;t) - \mathbf{1}_{t \Omega_{M^{-1}\gamma}}(x)|=| \int_{y \notin t \Omega_{\gamma}} \frac{t^{2}}{4 \pi} e^{- \frac{\lVert Mx - y \rVert ^{2}}{4}} dy| \textit{ if } Mx \in t \Omega_{\gamma}$$
and 
$$|\chi_{\gamma, M}(x;t) - \mathbf{1}_{t \Omega_{M^{-1}\gamma}}(x)|=| \int_{y \in t \Omega_{\gamma}} \frac{t^{2}}{4 \pi} e^{- \frac{\lVert Mx - y \rVert ^{2}}{4}} dy| \textit{ if } Mx \notin t \Omega_{\gamma}$$
because of Equation $(\ref{eq12})$ and by making the change of variable $y = Mu$. \\
The proof of Lemma 3.2 from $\cite{bleher1992distribution}$ gives the wanted result.
\end{proof}
The second sublemma gives an estimate of $\text{dist}(Mn,t \gamma)$. To state it, we need some notations. Like in $\cite{bleher1992distribution}$, let the curve $\gamma$ be defined in the polar coordinates $(r,\varphi)$ by the equation 
\begin{equation} 
\label{eq100001}
r= \Gamma(\varphi) \textit{.}
\end{equation}
Let's define : 
\begin{equation} 
\label{eq100002}
r_{\gamma}(x)= \frac{\lVert x \rVert}{\Gamma(\varphi(x))} 
\end{equation}
where $\varphi(x)$ is the angular coordinate of $x$.
Then, one has that there exists $C > 0$ small enough so that for every $x \in \mathbb{R}^{2}$, 
\begin{equation} 
\label{eq100003}
\text{dist}(x, t \gamma) \geqslant C |r_{\gamma}(x) - t| \textit{.}
\end{equation}
We deduce the following sublemma : 
\begin{sublemma}
\label{lemme4}
For all $L \in \mathscr{S}_{2}$, for all $t > 0$, for all $n \in \mathbb{Z}^{2}$, we have that : 
$$\text{dist}(Mn,t \gamma) \geqslant C | r_{\gamma}(M n) - t| \textit{.}$$
\end{sublemma}
Now we can prove Lemma $\ref{lemme32}$.
\begin{proof}[Proof of Lemma $\ref{lemme32}$]
By using the Equation ($\ref{eq19}$), we have that : 
$$\Delta_{1}(L,t) \leqslant \frac{1}{\sqrt{t}}\sum_{n \in \mathbb{Z}^{2}} |\chi_{\gamma, M}(n;t) - \mathbf{1}_{t \Omega_{M^{-1}\gamma}}(n)| $$ 
because, also, $\mathcal{R}(t \Omega_{\gamma},L)= \mathcal{R}(t \Omega_{M^{-1}\gamma},\mathbb{Z}^{2})$.
So, Sublemma $\ref{lemme3}$ and Sublemma $\ref{lemme4}$ imply that : 
$$\Delta_{1}(L,t) \leqslant \frac{1}{\sqrt{t}} \sum_{n \in \mathbb{Z}^{2}} e^{-C^{2} \frac{t^{2}}{4} | r_{\gamma}(M n) - t|^{2} } \textit{.} $$
The essential part of the right-hand side are the terms such that $n \in \mathbb{Z}^{2}$ verify that
\begin{equation}
\label{eq100004}
| r_{\gamma}(M n) - t| \leqslant \frac{1}{t^{3/4}} \textit{.} 
\end{equation}
Yet the number of $n \in \mathbb{Z}^{2}$ that belong to such an annulus is of order $t^{\frac{1}{4}}$. \\
So, one has finally : 
$$\Delta_{1}(L,t) = O (\frac{1}{t^{\frac{1}{4}}}) \textit{.} $$
\end{proof}
\subsubsection{Proof of Lemma $\ref{lemme33}$}
To prove Lemma $\ref{lemme33}$, we first need to give another expression of $F(M,t)$, obtained via the Poisson formula. It is the object of the following lemma : 
\begin{lemma}
\label{lemme34}
\begin{align}
& F(M,t) = \frac{1}{2\pi}\sum_{n \in \mathbb{Z}^{2}-\{ 0 \}}  \frac{ \widetilde{\lambda_{M}}(2 \pi n ; t)}{\lVert (M^{-1})^{T}n \rVert^{\frac{3}{2}}} \\
& \left(\rho_{\gamma}((M^{-1})^{T}n) \cos( 2 \pi t Y_{\gamma}((M^{-1})^{T}n) - \frac{3 \pi}{4}) + \rho_{\gamma}(-(M^{-1})^{T}n) \cos(2 \pi t Y_{\gamma}(-(M^{-1})^{T}n) - \frac{3 \pi}{4}) \right) \nonumber \\
&   + O_{M}(t^{-1})  \nonumber
\end{align}
which becomes in the symmetric case 
$$
 F(M,t) = \frac{1}{\pi}\sum_{n \in \mathbb{Z}^{2}-\{ 0 \}}  \frac{ \widetilde{\lambda_{M}}(2 \pi n ; t)}{\lVert (M^{-1})^{T}n \rVert^{\frac{3}{2}}} \rho_{\gamma}((M^{-1})^{T}n) \cos( 2 \pi t Y_{\gamma}((M^{-1})^{T}n) - \frac{3 \pi}{4})   + O_{M}(t^{-1})  
$$
where the $M$ in index of $O_{M}$ is to signal that it depends on $M$ (or, equivalently, on the lattice $L$). 
\end{lemma}
To prove it, we first need a calculatory sublemma : 
\begin{sublemma}
\label{lemme2}
Let $\gamma$ be a simple, closed, analytic, strictly convex curve such that $0 \in \Omega_{\gamma}$. Let $D \in SL_{2}(\mathbb{R})$. Let $\tilde{\gamma} = D \gamma$. Then one has for every $\xi \in \mathbb{R}^{2}-\{ 0 \}$ :
$$x_{\tilde{\gamma}}(\xi) = D x_{\gamma}(D^{T} \xi) \textit{, }$$
$$\rho_{\tilde{\gamma}}(\xi) = \frac{\lVert  \xi \rVert^{3}}{\lVert D^{T} \xi \rVert^{3}} \rho_{\gamma}(D^{T} \xi) \textit{ and } $$
$$ Y_{\tilde{\gamma}}(\xi) = Y_{\gamma}(D^{T} \xi) $$ 
where $D^{T}$ is the transpose of the matrix $D$.
\end{sublemma}
\begin{proof}
Let's call $T$ the unit tangent vector at $x_{\gamma}(D^{T} \xi)$ such that $(D^{T} \xi, T)$ is a orthogonal and direct basis of $\mathbb{R}^{2}$. Then, $D T$ is a unit tangent vector at $D x_{\gamma}(D^{T} \xi)$.\\
Let's call $k$ the unit normal exterior vector of $\tilde{\gamma}$ at $D x_{\gamma}(D^{T} \xi)$. Then, one knows that $ k$ is orthogonal to $D T$ and $(k,DT)$ is direct because $\det(D) = 1$. By property of the adjoint operator, one knows that $D^{T} k$ is orthogonal to $T$ and thus there exists $\alpha \in \mathbb{R}-\{ 0 \}$ such that 
\begin{equation}
\label{eq100005}
D^{T} k = \alpha D^{T} \xi \textit{.}
\end{equation}
So, one gets that : 
\begin{equation}
\label{eq100006}
k = \alpha \xi \textit{.}
\end{equation}
Yet, one has also that $(k,DT)$ and $(D^{T} \xi, T)$ are direct and orthogonal basis of $\mathbb{R}^{2}$. So, one must have $\alpha > 0$ and it gives us the first wanted result. \\
Now, concerning the third equality, one has that, by definition :  
\begin{equation}
\label{eq100007}
Y_{\tilde{\gamma}}(\xi) = <\xi, x_{\tilde{\gamma}}(\xi)>
\end{equation}
So, the first equality of Sublemma $\ref{lemme2}$ gives us that : 
\begin{equation}
\label{eq100008}
Y_{\tilde{\gamma}}(\xi) = <\xi, D x_{\gamma}(D^{T} \xi) >  \textit{.}
\end{equation}
By using the adjoint property, one finds the wanted result : 
\begin{equation}
\label{eq100009}
Y_{\tilde{\gamma}}(\xi) = <D^{T} \xi,  x_{\gamma}(D^{T} \xi) > =  Y_{\gamma}(D^{T} \xi) \textit{.}
\end{equation}
So, one gets the third equality. \\
Concerning the second equality, one knows that $ t \longmapsto \gamma(t)$ is a parametrization of the curve $\gamma$ and that $t \longmapsto D \gamma(t)$ is a parametrization of the curve $D \gamma$. So, one can use these parametrizations to compute $\rho_{\gamma}$ and $\rho_{\tilde{\gamma}}$. \\
By using the fact that $D \in SL_{2}(\mathbb{R})$, one has that : 
\begin{equation}
\label{eq1000010}
\rho_{\tilde{\gamma}}(D \gamma(t)) = (\frac{\lVert D \gamma'(t) \rVert}{\lVert \gamma'(t) \rVert})^{\frac{3}{2}} \rho_{\gamma}(\gamma(t)) \textit{.} 
\end{equation}
Let's call $t_{0}$ the instant such that $$\gamma(t_{0}) = x_{\gamma}(D^{T} \xi)$$ and so, according to the first result of the Sublemma $\ref{lemme2}$, one has that 
\begin{equation}
\label{eq1000011}
 D \gamma(t_{0}) = x_{\tilde{\gamma}}(\xi) \textit{.} 
\end{equation}
Let's set $\gamma(t) = (x(t),y(t))$ and $ D \gamma(t) = (\phi_{1}(t), \phi_{2}(t))$. \\
Then one has at the instant $t=t_{0}$ : 
\begin{equation}
\label{eq1000012}
\alpha R D^{T} \xi = (x'(t_{0}),y'(t_{0}))
\end{equation}
and
\begin{equation}
\label{eq1000013}
\beta R \xi = (\phi_{1}'(t_{0}),\phi_{2}'(t_{0}))
\end{equation}
with $\alpha, \beta > 0$ and $R = \begin{pmatrix} 0 & -1 \\ 1 & 0 \end{pmatrix} \in SO_{2}(\mathbb{R})$. \\
By using the fact that $D R D^{T} = R$, one finds that :
\begin{equation}
\label{eq1000014}
\alpha = \beta \textit{.}
\end{equation}
By using Equation ($\ref{eq1000011}$), Equation ($\ref{eq1000012}$), Equation ($\ref{eq1000013}$), Equation ($\ref{eq1000014}$) and by using Equation ($\ref{eq1000010}$) at the instant $t = t_{0}$, one has the third wanted equality.
\end{proof}
We can now tackle the proof of Lemma $\ref{lemme34}$.
\begin{proof}[Proof of Lemma $\ref{lemme34}$]
According to the Equation ($\ref{eq16}$), the Poisson summation formula and because of the fact that $\widetilde{\mathbf{1}_{t \Omega_{M^{-1}\gamma}}} (0) = \text{Area}(t \Omega_{M^{-1}\gamma})$ one has that : 
\begin{equation}
\label{eq26}
F(M,t) = \frac{1}{\sqrt{t}} \sum_{n \in \mathbb{Z}^{2}-\{ 0 \}} \widetilde{\mathbf{1}_{t \Omega_{M^{-1}\gamma}}}(2 \pi n) \widetilde{\lambda_{M}}(2 \pi n ; t) \textit{.}
\end{equation}
Yet, according to Lemma 2.1 from $\cite{bleher1992distribution}$, one has that :
\begin{equation}
\label{eq27}
  \widetilde{\mathbf{1}_{t \Omega_{M^{-1} \gamma}}}(\xi) = \sqrt{t} \lVert \xi \rVert^{-\frac{3}{2}} \sum_{\pm} \sqrt{2 \pi \rho_{M^{-1} \gamma}( \pm \xi)} \exp (\pm i( t Y_{ M^{-1} \gamma}(\pm \xi) - \frac{3 \pi}{4})) + O_{M}(t^{- \frac{1}{2}} \lVert \xi \rVert ^{-\frac{5}{2}})
\end{equation}
By using Sublemma $\ref{lemme2}$ with $M^{-1} = D$, we get with Equation ($\ref{eq26}$) and by grouping the $n$ and $-n$ terms in the Fourier series, we get the wanted result.
\end{proof}
Let's set 
\begin{equation}
\label{eq1000015}
\nu(l,t) = \rho_{\gamma}(l) \cos( 2 \pi t Y_{\gamma}(l) - \frac{3 \pi}{4}) \textit{.}
\end{equation}
Using the Equation ($\ref{eq13}$), the Equation ($\ref{eq19}$), the fact that if $M$ represents a lattice $L$, $(M^{-1})^{T}$ represents the dual lattice $L^{\perp}$, and the previous lemma, that is Lemma $\ref{lemme34}$, one gets that : 
\begin{equation}
\label{eq29}
(\Delta_{2})_{A}(L,t) \leqslant \Delta_{2,1}(L,t) + (\Delta_{2,2})_{A}(L,t) + (\Delta_{2,3})_{A}(L,t)
\end{equation}
where 
\begin{equation}
\label{eq30}
\Delta_{2,1}(L,t) = O_{M}(t^{-1})
\end{equation}
\begin{equation}
\label{eq31}
(\Delta_{2,2})_{A}(L,t) =\frac{1}{2 \pi} |  \sum_{\substack{ l \in L^{\perp} \\ 0 < \lVert l \rVert \leqslant A }} \frac{1}{\lVert l \rVert^{\frac{3}{2}}} (\nu(l,t) + \nu(-l,t) )(1 - e^{- (2 \pi)^{2} \frac{\lVert l \rVert^{2}}{t^{2}}}) |
\end{equation}
\begin{equation}
\label{eq32}
(\Delta_{2,3})_{A}(L,t) =\frac{1}{2 \pi}  |\sum_{\substack{ l \in L^{\perp} \\ A < \lVert l \rVert }} \frac{1}{\lVert l \rVert^{\frac{3}{2}}} (\nu(l,t) + \nu(-l,t) ) e^{- (2 \pi)^{2} \frac{\lVert l \rVert^{2}}{t^{2}}}| \textit{.}
\end{equation}
So, if we prove the following lemmas, we will get Lemma $\ref{lemme33}$ and, $\textit{in fine}$, get Lemma $ \ref{lemme30}$ : 
\begin{lemma} 
\label{lemme35}
$\Delta_{2,1}(L,t)$ converges almost surely to $0$ when $t \rightarrow \infty$.
\end{lemma}
Let's remark, by the way, that this last lemma is immediate according to equation $(\ref{eq30})$.
\begin{lemma}
\label{lemme36}
For all $A > 0$, $(\Delta_{2,2})_{A}(L,t)$ converges to $0$ when $t \rightarrow \infty$.
\end{lemma}
\begin{lemma}
\label{lemme37}
For all $\alpha > 0$, for all $A$ large enough, for all $t$ large enough,
$$\mathbb{P}((\Delta_{2,3})_{A}(L,t) \geqslant \alpha) \leqslant \alpha \textit{.}$$
\end{lemma}
\begin{proof}[Proof of Lemma $\ref{lemme33}$]
Let $\alpha > 0$. Let's take $A$ large enough so that for all $t$ large enough, $$\mathbb{P}((\Delta_{2,3})_{A}(L,t) \geqslant \alpha) \leqslant \alpha \textit{.}$$ 
It is possible according to Lemma $\ref{lemme37}$. \\
According to Lemma $\ref{lemme36}$, according to Lemma $\ref{lemme35}$ and because the almost-sure convergence imply the convergence in probability, even if it means taking $t$ larger, one can suppose that : $$\mathbb{P}((\Delta_{2,1})(L,t) \geqslant \alpha) \leqslant \alpha \textit{ and}$$ 
$$\mathbb{P}((\Delta_{2,2})_{A}(L,t) \geqslant \alpha) \leqslant \alpha \textit{.}$$ 
By using equation $\ref{eq29}$, one gets the wanted result.
\end{proof}
Before following with the proof of Lemma $\ref{lemme30}$, let's say a few words about Lemma $\ref{lemme36}$ and Lemma $\ref{lemme37}$. The first tells us that the non-regularized Fourier series is "close" enough to the regularized Fourier series whereas the second one tells us that the large terms of the regularized Fourier series do not matter, in a certain sense, for our study. \\
It remains only to prove Lemma $\ref{lemme36}$ and Lemma $\ref{lemme37}$. Because the density of $\tilde{\mu}_{2}$ is bounded, we only need to prove these lemmas for $\tilde{mu}_{2} = \mu_{2}$ and we will make this assumption for the rest of the section. We are now going to prove Lemma $\ref{lemme36}$ and Lemma $\ref{lemme37}$. 
\subsubsection{Proof of Lemma $\ref{lemme36}$}
\begin{proof}[Proof of Lemma $\ref{lemme36}$]
Let $l \in L^{\perp}$. Then one has : 
\begin{equation}
\label{eq33}
|1 - e^{- (2 \pi)^{2} \frac{\lVert l \rVert^{2}}{t^{2}}} | \leqslant \frac{(2 \pi)^{2} \lVert l \rVert^{2}}{t^{2}} \textit{.}
\end{equation}
With this equation and with Equation $(\ref{eq31})$, one gets that : 
\begin{equation}
\label{eq34}
(\Delta_{2,2})_{A}(L,t) \leqslant \sum_{\substack{l \in L^{\perp} \\ 0 < \lVert l \rVert \leqslant A}} \frac{M}{t^{2}} \lVert l \rVert^{\frac{1}{2}} 
\end{equation}
with $M > 0$ because $\nu$ is a bounded function. 
It follows that there exists $C(L) > 0$ such that : 
\begin{equation}
\label{eq35}
(\Delta_{2,2})_{A}(L,t) \leqslant C(L) \frac{A^{\frac{5}{2}}}{t^{2}} \textit{.}
\end{equation}
\end{proof}
\subsubsection{Proof of Lemma $\ref{lemme37}$} 
To prove Lemma $\ref{lemme37}$ we need to use what are called $\textit{Siegel}$ and $\textit{Rogers formulas}$. Theses formulas will also be useful later in this paper. \\
By setting $c_{k} = \zeta(2)^{-k}$ for $k$ an integer larger than $1$ and where $\zeta$ denotes the $\zeta$ function of Riemann, one has the following formulas :
\begin{lemma}[\cite{marklof1998n},\cite{vinogradov2010limiting},\cite{kelmer2021second}]
\label{lemme5} 
For $f$ a piecewise smooth function with compact support on $\mathbb{R}^{2}$, one has : \begin{itemize}
\item $$  \int_{\mathscr{S}_{2}} \mathcal{S}(f) d \mu_{2} = c_{1} \int_{\mathbb{R}^{2}} f d\lambda  $$
\item When $f$ is even, 
\begin{align*}
(b) \int_{\mathscr{S}_{2}} \mathcal{S}(f)^{2} d \mu_{2}  \leqslant C \int_{\mathbb{R}^{2}} f^{2} d\lambda + c_{2}(\int_{\mathbb{R}^{2}}f d\lambda)^{2} 
\end{align*}
where $C > 0$. 
\end{itemize}

\end{lemma}
With this lemma, we are going to prove two lemmas that will enable us to prove Lemma $\ref{lemme37}$ by using Chebyshev's inequality : the first one is intended to estimate the expectation of $(\Delta_{2,3})_{A}(L,t)$ to see that it goes to $0$ when $t \rightarrow \infty$ (uniformly in $A$), the second one is intended to estimate its variance to see that it can be as uniformly small in $t$ as one wants if $A$ is chosen large enough. Until the end of this section, we are going to suppose $A > 1$.
\begin{lemma}
\label{lemme6}
$$\mathbb{E}((\Delta_{2,3})_{A}(L,t)) = O(\frac{1}{t }) \textit{.}$$
\end{lemma}
\begin{proof}
One has :
\begin{equation}
\label{eq36}
(\Delta_{2,3})_{A}(L,t) =  \sum_{\substack{ l \in L^{\perp} \\ A < \lVert l \rVert }} f(l) 
\end{equation}
where 
\begin{equation}
\label{eq37}
f(l) = \frac{1}{2 \pi}\frac{1}{\lVert l \rVert^{\frac{3}{2}}}  (\nu(l,t) + \nu(-l,t) ) e^{- (2 \pi)^{2} \frac{\lVert l \rVert^{2}}{t^{2}}} \textit{.}
\end{equation}
The Lemma $\ref{lemme5}$ gives us then that : 
\begin{equation}
\label{eq38}
| \mathbb{E}((\Delta_{2,3})_{A}(L,t)) |= C | \int_{\mathbb{R}^{2}} f(x) \mathbf{1}_{ \lVert x \rVert > A } dx | \textit{.}
\end{equation}
 \\ 
By passing into polar coordinates $(r,\theta)$, one gets that : 
\begin{equation}
\label{eq39}
| \mathbb{E}((\Delta_{2,3})_{A}(L,t)) | = 2 C | \int_{\theta = 0}^{2 \pi} \int_{r \geqslant A}   \frac{ \cos(2 \pi t r h(\theta) -  \frac{3 \pi}{4})  e^{- (2  \pi)^{2} \frac{r}{t^{2}}}}{ r^{\frac{1}{2}}}  dr d \theta |
\end{equation}
by setting $h(\theta) = Y_{\gamma}((\cos(\theta),\sin(\theta))$ and by using the fact that $Y_{\gamma}$ is positively homogeneous. 
Furthermore, an integration by part gives us that : 
\begin{align}
\label{eq41}
\int_{r \geqslant A}   \frac{ \cos(2 \pi t r h(\theta) -  \frac{3 \pi}{4})  e^{- (2  \pi)^{2} \frac{r}{t^{2}}}}{ r^{\frac{1}{2}}}  dr & \\
 = - \frac{e^{- \frac{(2 \pi)^{2} r^{2} }{t^{2}}}  \sin(2 \pi t A h(\theta) - \frac{3 \pi}{4})}{2 \pi t h(\theta) A^{\frac{1}{2}}} & \nonumber  \\
 + \frac{1}{4 \pi t h(\theta) } \int_{r \geqslant A} \sin(2 \pi t r h(\theta)- \frac{3 \pi}{4}) \frac{e^{- \frac{(2 \pi)^{2} r^{2} }{t^{2}}} }{r^{\frac{3}{2}}} dr & \nonumber \\
 + \frac{4 \pi}{t^{3} h(\theta) } \int_{r \geqslant A}   \sin(2 \pi t r h(\theta)- \frac{3 \pi}{4}) e^{- \frac{(2 \pi)^{2} r^{2}}{t^{2}}} r^{\frac{1}{2}} dr & \nonumber \textit{.}
\end{align}
By using that $r^{\frac{1}{2}} \leqslant r$ (because $A > 1$) and the fact that $h(\theta)$ admits a positive lower bound (because $0 \in \Omega_{\gamma}$), by estimating the three terms of the right member, one gets that : 
\begin{equation}
\label{eq4100}
| \int_{r \geqslant A}   \frac{ \cos(2 \pi t r h(\theta) -  \frac{3 \pi}{4})  e^{- (2  \pi)^{2} \frac{r}{t^{2}}}}{ r^{\frac{1}{2}}}  dr | \leqslant \frac{C}{t}
\end{equation}
with $C > 0$ that does not depend on $\theta$. \\
By using Equation ($\ref{eq39}$) and Equation ($\ref{eq4100}$), one gets that : 
$$\mathbb{E}((\Delta_{2,3})_{A}(L,t))= O(\frac{1}{t}) \textit{.}$$
\end{proof}
\begin{lemma}
\label{lemme7}
$$\text{Var}((\Delta_{2,3})_{A}(L,t)) = O(\frac{1}{A})  $$
where the $O$ can be chosen independent from $t$. 
\end{lemma}
\begin{proof}
By using the same notation as before, by using again the Lemma $\ref{lemme5}$ and by using the Lemma $\ref{lemme6}$, one gets that : 
\begin{equation}
\label{eq42}
\text{Var}((\Delta_{2,3})_{A}(L,t)) \leqslant C \int_{\mathbb{R}^{2}} f^{2}(x) \mathbf{1}_{\lVert x \rVert > A} dx  \textit{.}
\end{equation}
So, by passing into polar coordinates and by using the fact that $\nu$ is bounded, one gets that : 
\begin{equation}
\label{eq43}
\text{Var}((\Delta_{2,3})_{A}(L,t))\leqslant C \int_{r=A}^{\infty} \frac{1}{r^{2}}  dr  \textit{.}
\end{equation}
By integrating, we get the wanted result. 
\end{proof}
We can now prove the Lemma $\ref{lemme37}$.
\begin{proof}[Proof of Lemma $\ref{lemme37}$]
The Chebyshev's inequality gives the wanted result if, first, we choose $A$ large enough and, second, we choose $t$ large enough so that $\mathbb{E}((\Delta_{2,3})_{A}(L,t))$ and $\text{Var}((\Delta_{2,3})_{A}(L,t))$ are small enough. These choices are possible according to Lemmas $\ref{lemme6}$ and $\ref{lemme7}$.
\end{proof}
So, now the proof of Lemma $\ref{lemme30}$ is complete and we will conclude this section by proving the Lemma $\ref{lemme31}$ so that the proof of Proposition $\ref{prop8}$ will be complete.
\subsection{Proof of Lemma $\ref{lemme31}$}
To prove the Lemma $\ref{lemme31}$, we are going to take the same kind of approach as before : estimate the expectation and the variance of the quantity $S_{A}-S_{A,prime}$ and get the result via the Chebyshev's inequality. \\
We have that : 
\begin{equation}
\label{eq45}
H_{A}(L,t)-S_{A,prime}(L,t) = \sum_{l \in L^{\perp} \textit{prime}} f(l) 
\end{equation}
where
\begin{equation}
\label{eq46}
f(l) = \frac{1}{2 \pi \lVert l \rVert^{\frac{3}{2}}} \mathbf{1}_{0 < \lVert l \rVert \leqslant A} \sum_{k \geqslant \lfloor \frac{A}{\lVert l \rVert} \rfloor + 1} \frac{(\nu(k l,t) + \nu(-k l,t) )}{k^{\frac{3}{2}}} \textit{.}
\end{equation}
With this expression, we see that we are going to have a little problem of integrability at $0$ if we use the Lemma $\ref{lemme5}$. That's why, we have to exclude $0$ and we will suppose that $L$ is chosen so that $\lVert L^{\perp} \rVert_{1} \geqslant \epsilon$ where $0 < \epsilon < 1$. Only a small number of lattices is excluded according to this lemma : 
\begin{lemma}
\label{lemme8}
For every $0 < \epsilon < 1$, one has that 
$$\mathbb{P}(\lVert L \rVert_{1} < \epsilon) = O(\epsilon^{2}) \textit{.} $$
\end{lemma}
\begin{proof}
It is a consequence of Lemma $\ref{lemme5}$ by taking 
$$\mathcal{S}(f)(L) = \sum_{l \in L} \mathbf{1}_{B_{f}(0,\epsilon)}(l) $$ 
where $ \mathbf{1}_{B_{f}(0,\epsilon)}(l)$ is the indicator function of the closed ball for the norm $\lVert \cdot \rVert$ centred on $0$ of radius $\epsilon$.
\end{proof}
Thus, for the chosen lattices, we have :
\begin{equation}
\label{eq47}
H_{A}(L,t)-S_{A,prime}(L,t) = \sum_{l \in L^{\perp} \textit{prime}} f(l) \mathbf{1}_{\lVert l \rVert \geqslant \epsilon} = \Delta_{3,\epsilon,A,t}(L) 
\end{equation}
(this equation defines $\Delta_{3,\epsilon,A,t}(L) $).
\begin{lemma}
\label{lemme9}
$$\mathbb{E}(\Delta_{3,\epsilon,A,t}(L)) = O_{\epsilon,A}(\frac{1}{t}) \textit{.}$$
\end{lemma}
\begin{proof}
By using the Lemma $\ref{lemme5}$, one gets that : 
\begin{equation}
\label{eq48}
\mathbb{E}(\Delta_{3,\epsilon,A,t}(L)) \leqslant C \int_{r=\epsilon}^{A} \int_{\theta = 0 }^{2 \pi} \frac{1}{r^{\frac{1}{2}}} \sum_{k \geqslant \frac{A}{r}} \frac{\nu( k r (\cos(\theta), \sin(\theta)),t) }{k^{\frac{3}{2}}} dr \textit{.} 
\end{equation}
Lebesgue's dominated convergence theorem gives us that : 
\begin{equation}
\label{eq49}
\int_{\epsilon}^{A} \frac{1}{r^{\frac{1}{2}}} \sum_{k \geqslant \frac{A}{r}} \frac{\nu( k r (\cos(\theta), \sin(\theta)),t)}{k^{\frac{3}{2}}} dr = \sum_{k=1}^{\infty} \frac{1}{k^{\frac{3}{2}}} \int_{\max(\frac{A}{k}, \epsilon)}^{A} \frac{\nu( k r (\cos(\theta), \sin(\theta)),t)}{r^{\frac{1}{2}}}dr \textit{.}
\end{equation}
An integration by part (on the variable $r$) as in the proof of Lemma $\ref{lemme6}$ and the Equation $(\ref{eq48})$ give us finally that :  
\begin{equation}
\label{eq50}
\mathbb{E}(\Delta_{3,\epsilon,A,t}(L)) = O_{\epsilon,A}(\frac{1}{t}) \textit{.}
\end{equation}
\end{proof}
\begin{lemma}
\label{lemme10}
There exists $K > 0$ such that : 
$$\text{Var}(\Delta_{3,\epsilon,A,t}(L)) \leqslant K(- \frac{\log(\epsilon)}{A} + \frac{\log(A)}{A}) \textit{.} $$
\end{lemma}
\begin{proof}
Lemma $\ref{lemme5}$ gives us that : 
\begin{equation}
\label{eq51}
\text{Var}(\Delta_{3,\epsilon,A,t}(L)) \leqslant C  \int_{r=\epsilon}^{A} \int_{\theta = 0}^{2 \pi} \frac{1}{r^{2}} (\sum_{k \geqslant \frac{A}{r}} \frac{\nu( k r (\cos(\theta), \sin(\theta)),t)}{k^{\frac{3}{2}}})^{2} dr \textit{.}
\end{equation}
Yet one also has that for all $x > 0$ : 
\begin{equation}
\label{eq52}
\sum_{k \geqslant x} \frac{1}{k^{\frac{3}{2}}} \leqslant \frac{D}{x^{\frac{1}{2}}}
\end{equation}
where $D >0$. \\
Thus, Equation $(\ref{eq51})$, Equation $(\ref{eq52})$ and the fact that $\nu$ is bounded imply that : 
\begin{equation}
\label{eq53}
\text{Var}(\Delta_{3,\epsilon,A,t}(L)) \leqslant \frac{2 \pi C D}{A} \int_{\epsilon}^{A} \frac{1}{r} dr \textit{.}
\end{equation}
\end{proof}
We can now give the proof of Lemma $\ref{lemme31}$.
\begin{proof}[Proof of Lemma $\ref{lemme31}$]
First we take $1> \epsilon > 0$ small enough so that the measure of the neglected lattices, $\textit{id est}$ the lattices such that $\lVert L \rVert_{1} < \epsilon$, is small enough. It is possible according to Lemma $\ref{lemme8}$. \\
Then we take $A$ large enough so that $\text{Var}(\Delta_{3,\epsilon,A,t}(L))$ is small enough. It is possible according to Lemma $\ref{lemme10}$. \\
Finally, we take $t$ large enough so that $\mathbb{E}(\Delta_{3,\epsilon,A,t}(L))$ is small enough, which is possible according to Lemma $\ref{lemme9}$, and conclude by using Chebyshev's inequality.
\end{proof}
So, we are now brought back to the study of $S_{A,prime}(L,t)$ when $t \rightarrow \infty$ and the next section is dedicated to it. \\
$\textbf{We are going to replace}$ $L^{\perp}$ $\textbf{by}$ $L$ (it changes nothing because we are studying the asymptotic convergence in law with $L \in \mathscr{S}_{2}$ distributed according to $\tilde{\mu}_{2}$).
\section{Study of $S_{A,prime}(L,t)$ when $t \rightarrow \infty$}
\subsection{Reductions for the study of $S_{A,prime}(L,t)$ and proof of Theorem $\ref{thm2}$}
Before entering in the main object of this section, we need to do a small rewriting of $S_{A,prime}(L,t)$.\\
We recall that a vector $l \in L$ is prime if, and only if, $K l \in K L$ is prime where $K \in SL_{2}(\mathbb{R})$. Furthermore, a vector $(l_{1},l_{2}) \in \mathbb{Z}^{2}$ is prime if, and only if, $l_{1} \wedge l_{2} = 1$. \\
By using the symmetry $l \longmapsto -l$, we deduce that $S_{A,prime}(L,t)$ can be rewritten as followed : 
\begin{equation}
\label{eq55}
S_{A,prime}(L,t) = \frac{1}{\pi} \sum_{k \in \Pi_{A}(L)} \frac{Z_{k}(L,t)}{W_{k}(L)} 
\end{equation} 
where, for $k=(k_{1},k_{2}) \in \mathbb{Z}^{2}$,
\begin{equation}
\label{eq56}
W_{k}(L) = \lVert k_{1} e_{1}(L) + k_{2} e_{2}(L) \rVert^{\frac{3}{2}} \textit{, }
\end{equation}
\begin{equation}
\label{eq57}
Z_{k}(L,t) = \sum_{m \in \mathbb{N}-\{ 0 \}} \frac{  \nu( m (k_{1} e_{1}(L) + k_{2} e_{2}(L)),t) + \nu( -m (k_{1} e_{1}(L) + k_{2} e_{2}(L)),t)  }{m^{\frac{3}{2}}} \textit{, }
\end{equation}
and where
\begin{equation}
\label{eq58}
\Pi_{A}(L)= \{ (k_{1},k_{2}) \in \mathbb{Z}^{2} \text{ } | \text{ } k_{1} \wedge k_{2} = 1 \textit{, } k_{1} \geqslant 0 \textit{ } \lVert k_{1} e_{1}(L) + k_{2} e_{2} (L) \rVert \leqslant A \} 
\end{equation}
and here we agree that if $(k_{1},k_{2}) \in \Pi_{A}(L)$ then $k_{1} = 0$ implies that $k_{2} =1$ (for the definition of $e_{1}(L)$ and $e_{2}(L)$ see Definition $\ref{def1}$). \\
\\
Let's recall that : 
\begin{equation}
\label{eq67}
\Pi = \{ (k_{1},k_{2}) \in \mathbb{Z}^{2} \text{ } | \text{ } k_{1} \wedge k_{2} = 1 \textit{, } k_{1} \geqslant 0  \} \textit{.}
\end{equation}
Our goal now is to prove the following proposition : 
\begin{prop}
\label{prop10}
In the symmetric case, $ \{  Z_{k}(L,t) \}_{k \in \Pi}$ converge, when $t \rightarrow \infty$, in distribution towards $ \{ \rho_{\gamma}(k_{1} e_{1}(L) + k_{2} e_{2}(L)) \tilde{Z}_{k}(\omega) \}_{k \in \Pi}$ where $\tilde{Z}_{k}(\omega)$, with $\omega \in \Omega$, are independent identically distributed real random variables that have a compact support, are symmetrical and are non-zero. \\
In the non symmetric case, $ \{  Z_{k}(L,t) \}_{k \in \Pi}$ converge, when $t \rightarrow \infty$, in distribution towards $ \{ \rho_{\gamma}(k_{1} e_{1}(L) + k_{2} e_{2}(L)) \tilde{Z}_{k}(\omega_{1}) + \rho_{\gamma}(-(k_{1} e_{1}(L) + k_{2} e_{2}(L))) \tilde{Z}_{k}(\omega_{2}) \}_{k \in \Pi}$ where $(\omega_{1}, \omega_{2}) \in \Omega \times \Omega $, where $\tilde{Z}_{k}(\omega)$ are independent identically distributed real random variables that have a compact support, are symmetrical and are non-zero. 
\end{prop}
In the next section we are going to consider, in the symmetric case, the sums of the type $$\tilde{S}_{A}(\omega,L) = \sum_{k \in \Pi_{A}(L)} \frac{\rho_{\gamma}(k_{1} e_{1}(L) + k_{2} e_{2}(L)) \tilde{Z}_{k}(\omega)}{\lVert k_{1} e_{1}(L) + k_{2} e_{2}(L) \rVert^{\frac{3}{2}}}   $$
and the sums of the type, in the non symmetric case,
$$\tilde{S}_{A}(\omega_{1},\omega_{2},L) = \sum_{k \in \Pi_{A}(L)} \frac{\rho_{\gamma}(k_{1} e_{1}(L) + k_{2} e_{2}(L)) \tilde{Z}_{k}(\omega_{1} ) + \rho_{\gamma}(-(k_{1} e_{1}(L) + k_{2} e_{2}(L))) \tilde{Z}_{k}(\omega_{2} )}{\lVert k_{1} e_{1}(L) + k_{2} e_{2}(L) \rVert^{\frac{3}{2}}}   $$
where $Z_{k}$ are non-zero real independent identically distributed random variables from $\Omega \ni \omega$ that are symmetrical and have a compact support and where $(\omega_{1}, \omega_{2}) \in \Omega \times \Omega$. Proposition $\ref{prop15}$ tells us that : 
\begin{prop}
\label{prop15}
The sums of these types : 
 \begin{itemize}
\item converge almost surely 
\item their respective limits are symmetrical and their expectations are equal to $0$
\item their respective limits admit moment of order $1+\kappa$ for any $0 \leqslant \kappa < \frac{1}{3}$ 
\item their respective limits do not admit a moment of order $\frac{4}{3}$ when $\sigma(L) \geqslant m$ where $m > 0$ and where $L$ belongs to an event of the form $( \lVert L \rVert < \alpha )$ with $\alpha > 0$. 
\item when there exists $\alpha > 0 $ such that $\tilde{\mu}_{2}(\{ L \in \mathscr{S}_{2} \text{ } | \text{ } \lVert L \rVert_{1} < \alpha \}) = 0$ then their respective limits admit moments of all order $1 \leqslant p < \infty$.  
\end{itemize}
\end{prop}
We are going to see now that it is enough to prove Proposition $\ref{prop10}$ and Proposition $\ref{prop15}$ to establish Theorem $\ref{thm2}$, with the exception of the exact form of the limiting law (yet it is given by Proposition $\ref{prop11}$).
\begin{proof}[Proof of Theorem $\ref{thm2}$]
Let $\psi \in C_{c}^{\infty}(\mathbb{R})$. Let $\epsilon > 0$. According to Proposition $\ref{prop8}$, we can take $A$ as large as we want and then $t$ as large as we want so that : 
\begin{equation}
\label{eq502}
| \mathbb{E}\left( \psi(\frac{\mathcal{R}(t \Omega_{\gamma}, L)}{\sqrt{t}}) \right) - \mathbb{E}\left( \psi(S_{A,prime}(L,t)) \right) | \leqslant \epsilon \textit{.} 
\end{equation}
Thanks to Proposition $\ref{prop10}$, one has also that : 
\begin{equation}
\label{eq503}
| \mathbb{E}\left( \psi(S_{A,prime}(L,t)) \right) - \mathbb{E}\left( \psi(\tilde{S}_{A}(\omega,L)) \right) | \leqslant \epsilon 
\end{equation}
where the $Z_{k}(\omega)$ in $\tilde{S}_{A}(\omega,L)$ are given by Proposition $\ref{prop10}$. \\
Furthermore, Proposition $\ref{prop15}$ gives us that : 
\begin{equation}
\label{eq504}
| \mathbb{E}\left( \psi(\tilde{S}_{A}(\omega,L)) \right) - \mathbb{E}\left( \psi(\lim_{A \rightarrow \infty} \tilde{S}_{A}(\omega,L)) \right) | \leqslant \epsilon 
\end{equation}
with $\lim_{A \rightarrow \infty} \tilde{S}_{A}(\omega,L)$ that verify all the listed properties. \\
So, Equation $(\ref{eq502})$, Equation $(\ref{eq503})$ and Equation $(\ref{eq504})$ give the wanted result.
\end{proof}
The main reason why the $Z_{k}$ are going to be independent from $L$ is the presence of the factor $t$. \\
The main reasons why the rest of Proposition $\ref{prop10}$ will be true are the presence of the factor $t$ in $Z_{k}$ and the fact that the coefficients of order 1 of the Taylor series of $(Y( k_{1} e_{1}(L) + k_{2} e_{2}(L)))_{k \in \Pi_{A}(L)}$ on a small geodesic segment are $\mathbb{Z}$-free and, in the non symmetric case, $(Y( k_{1} e_{1}(L) + k_{2} e_{2}(L)), Y(-(k_{1} e_{1}(L) + k_{2} e_{2}(L)))_{k \in \Pi_{A}(L)}$ on a small geodesic segment are $\mathbb{Z}$-free. \\
\\
In order to prove Proposition $\ref{prop10}$, it is actually enough to prove the following proposition by using the definition of $\nu$ (see Equation ($\ref{eq1000015}$)) : 
\begin{prop}
\label{prop11}
For $k=(k_{1},k_{2})$, let 
\begin{equation}
\label{eq59}
\theta_{k}(L,t) = t  Y_{\gamma}( k_{1} e_{1}(L) + k_{2} e_{2}(L))) \mod 1 \textit{.}
\end{equation}
Then, we have that
$\{ \theta_{k}(L,t) \}_{k \in \Pi}$ converge, when $t \rightarrow \infty$, towards random variable that are independent identically distributed, are distributed according to the Lebesgue measure $\lambda$ over $\mathbb{R}/\mathbb{Z}$ and are independent from $L$. \\
In the non symmetric case, we have that $\{ \theta_{k}(L,t) \}_{k \in \Pi}$ and $\{ \theta_{-k}(L,t) \}_{k \in \Pi}$ converge, when $t \rightarrow \infty$, towards random variable that are independent identically distributed, are distributed according to the Lebesgue measure $\lambda$ over $\mathbb{R}/\mathbb{Z}$ and are independent from $L$.
\end{prop}
$\textbf{Thanks to this proposition, we now understand why the limit law of }$ $ \frac{\mathcal{R}(t \mathbb{D}^{2},L)}{\sqrt{t}}$ $ \textbf{ is} $ \\
$\textbf{given by } S(\theta,L) \textit{.}$
To prove this last proposition, it is sufficient to prove the following proposition where $e(\theta)$ stands for $\exp(i 2 \pi \theta)$ : 
\begin{prop}
\label{prop12}
For every $l \in \mathbb{N}-\{0\}$, for every $\psi \in C_{c}^{\infty}(\mathscr{S}_{2})$, for every $(p_{1},\cdots,p_{l}) \in \mathbb{Z}^{l}-\{ 0 \}$, one has : 
\begin{equation}
\label{eq66}
\mathbb{E} \left( \psi(L) e(\sum_{h=1}^{l} p_{h} \theta_{k_{h}}) \right) \underset{t \rightarrow \infty}{\rightarrow} 0
\end{equation}
where the $k_{h} \in \Pi$ are all distinct. \\
In the non symmetric case, one has that for every $l \in \mathbb{N}-\{0\}$, for every $\psi \in C_{c}^{\infty}(\mathscr{S}_{2})$, for every $(p_{-l},\cdots,p_{-1},p_{1},\cdots,p_{l}) \in \mathbb{Z}^{2l}-\{ 0 \}$, one has : 
\begin{equation}
\label{eq66}
\mathbb{E} \left( \psi(L) e(\sum_{h=1}^{l} p_{h} \theta_{k_{h}} + \sum_{h=1}^{l} p_{-h} \theta_{-k_{h}} ) \right) \underset{t \rightarrow \infty}{\rightarrow} 0
\end{equation}
where the $k_{h} \in \Pi$ are all distinct. 
\end{prop}
Before passing to the proof of Proposition $\ref{prop12}$, let's give some heuristic about it in the symmetric case, the non symmetric case being similar here.\\ 
Basically, by working with a foliation of the space $\mathscr{S}_{2}$ given by small enough geodesic segments, we are first going to have : $$\mathbb{E} \left(  \psi(L) e(\sum_{h=1}^{l} p_{h} \theta_{k_{h}}) \right) \approx \mathbb{E}(\psi) \mathbb{E} \left(e(\sum_{h=1}^{l} p_{h} \theta_{k_{h}}) \right)  $$
due to the presence of the factor $t$ in $\theta$. \\
The right member will go to $0$ when $t$ goes to infinity because a Riemann-Lebesgue lemma will apply because quantities "close" to the variables $\theta_{k}$ are typically $\mathbb{Z}$-free (see the heuristic explanation of the second step). \\ 
The rest of this section is now dedicated to the proof of the Proposition $\ref{prop12}$.
\subsection{Foliation and local estimates}
We recall that a foliation of the space $\mathscr{S}_{2}$ is given by the orbits of the group $\delta$ where $$ \delta(\lambda) = \begin{pmatrix} \lambda & 0 \\ 0 & \frac{1}{\lambda} \end{pmatrix} \textit{.}$$
To prove Proposition $\ref{prop12}$, we are going to look at what it is happening on a small "segment" of the form 
\begin{equation}
\label{eq60}
J_{\epsilon}(L) = \{ \delta(\lambda) L \textit{ } | \textit{ } \lambda \in [\frac{1}{1+ \epsilon},1+\epsilon] \}
\end{equation}
where $L \in \mathscr{S}_{2}$ and $\epsilon > 0$ can be taken as small as possible. More precisely, we are going to show, when $t \rightarrow \infty$, the independence of the $(\theta_{k})$ and of $L$ over smalls segments of the form $J_{\epsilon}(L)$, as well as the fact that the $(\theta_{k})$ are identically distributed and distributed according to the normalized Haar measure over $\mathbb{R}/ \mathbb{Z}$.  \\
Let's call $s : (x,y) \in \mathbb{R}^{2} \longmapsto (x,-y)$ and let's set for $k=(k_{1},k_{2}) \in \Pi$ or $k \in - \Pi$ : 
\begin{equation}
\label{eq400}
W_{k}(L) = dY_{\gamma}(k_{1} e_{1}(L) + k_{2} e_{2}(L))\left( s(k_{1} e_{1}(L) + k_{2} e_{2}(L)) \right) 
\end{equation}
where $dY_{\gamma}(k_{1} e_{1}(L) + k_{2} e_{2}(L))(\cdot)$ stands for the differential of $Y_{\gamma}$ at the point $k_{1} e_{1}(L) + k_{2} e_{2}(L)$. \\
On a segment of the form $J_{\epsilon}(L)$, the following lemma basically tells us how we can estimate the quantities $Y_{\gamma}(k_{1}e_{1}(L) + k_{2} e_{2}(L))$ : 
\begin{lemma}
\label{lemme11}
For a typical $L \in \mathscr{S}_{2}$, there exists $\epsilon > 0$ small enough such that for every $\lambda \in [\frac{1}{1+\epsilon},1+\epsilon]$, $$e_{1}(\delta(\lambda)L) = \delta(\lambda)e_{1}(L) \textit{ and } e_{2}(\delta(\lambda)L) = \delta(\lambda)e_{2}(L) \textit{.}$$
Furthermore, for such a lattice $L$, for such $\lambda$, for $k=(k_{1},k_{2}) \in \Pi$ or for $k \in - \Pi$, we have for $h=\lambda-1$,
\begin{align}
\label{eq61}
& Y_{\gamma}(k_{1}e_{1}(\delta(\lambda)L) + k_{2} e_{2}(\delta(\lambda)L)) = Y_{\gamma}(k_{1} e_{1}(L) + k_{2} e_{2}(L) ) \\
& + h W_{k}(L)  + O_{k_{1},k_{2},L}(h^{2})   \nonumber \textit{.}
\end{align}
\end{lemma}
\begin{proof}
The first fact was proven in $\cite{trevisan2021limit}$ (see Lemma 15). \\
Let's note that, as $\gamma$ is analytical, $Y_{\gamma}$ is regular. As a consequence, the second fact is obtained from the first fact of Lemma $\ref{lemme11}$ and by a simple calculus of Taylor series. 
\end{proof}
To prove Proposition $\ref{prop12}$ we see, in light of Lemma $\ref{lemme11}$, that it would be convenient to prove the following proposition : 
\begin{prop}
\label{prop13}
For a typical $L \in \mathscr{S}_{2}$, for every $m \in \mathbb{N}-\{ 0 \}$, for every family $(p_{1},\cdots,p_{m}) \in \mathbb{Z}^{m}$, for every $k_{1},\cdots,k_{m} \in \Pi$ all distinct if 
\begin{equation}
\label{eq505}
\sum_{i=1}^{m}  p_{i} W_{k_{i}}(L) = 0
\end{equation}
then $p_{1} = \cdots = p_{m} = 0$. \\
In other words, for a typical $L \in \mathscr{S}_{2}$, $$\left( W_{k}(L) \right)_{k \in \Pi}$$ is a $\mathbb{Z}$-free family. \\
In the non symmetric case, for a typical $L \in \mathscr{S}_{2}$, for every $m \in \mathbb{N}-\{ 0 \}$, for every family $(p_{-m},\cdots,p_{-1},p_{1},\cdots,p_{m}) \in \mathbb{Z}^{2m}$, for every $k_{1},\cdots,k_{m} \in \Pi$ all distinct if 
\begin{equation}
\label{eq100010} 
\sum_{i=1}^{m}  p_{i} W_{k_{i}}(L) + \sum_{i=1}^{m}  p_{-i} W_{-k_{i}}(L)  = 0
\end{equation}
then $p_{-m}=\cdots=p_{-1}=p_{1} = \cdots = p_{m} = 0$. \\
In other words, in the non symmetric case, for a typical $L \in \mathscr{S}_{2}$, $$\left( W_{k}(L) \right)_{k \in -\Pi \cupdot \Pi}$$ is a $\mathbb{Z}$-free family. \\
\end{prop}
The next subsection is dedicated to prove this proposition. 
\subsection{Proof of Proposition $\ref{prop13}$}
To prove Proposition $\ref{prop13}$, we are following closely what was done in the Section 5 of $\cite{dolgopyat2014deviations}$ and we need four preliminary lemmas. \\
To state these lemmas we need to put in place some notations. Let's call 
\begin{equation}
\label{eq1000021}
P : X \in \mathbb{R}^{2} \longmapsto dY_{\gamma}(X)(s(X)) \textit{.}
\end{equation}
Let's also set for every $k=(k_{1},k_{2}) \in -\Pi \cupdot \Pi$, 
\begin{equation}
\label{eq1005}
f_{k} : (A=(X_{1},Y_{1}),B=(X_{2},Y_{2})) \in \mathbb{R}^{2} \times \mathbb{R}^{2} \longmapsto  P(k_{1} X + k_{2} Y)
\end{equation}
where $ X= (X_{1},X_{2})$ and $Y=(Y_{1},Y_{2}) \in \mathbb{R}^{2}$. \\
Let's set for every $L \in GL_{2}(\mathbb{R})$, 
\begin{equation} 
\label{eq1000016}
g_{L} : \delta \in \mathbb{R} \longmapsto P(L(1,\delta)) \textit{.}
\end{equation}
Let's also set for every $L \in GL_{2}(\mathbb{R})$, 
\begin{equation} 
\label{eq1000018}
\tilde{g}_{L} : \delta \in \mathbb{R} \longmapsto P(L(-1,-\delta)) \textit{.}
\end{equation}
\begin{lemma}
\label{lemme102}
If $\gamma$ is analytic, we have that for any $L \in GL_{2}(\mathbb{R})$, $g_{L}$ is (real) analytic and not equal to a polynomial.
\end{lemma}
\begin{proof}
First, let's note that we $X \in \mathbb{R}^{2} \longmapsto dY_{\gamma}(X)(s(X))$ is positively homogeneous because $X \longmapsto Y_{\gamma}(X)$ also is positively homogeneous. \\
So, one has that for every $\delta \in \mathbb{R}$, 
\begin{equation}
\label{eq1000019}
g_{L}(\delta) = \sqrt{1+\delta^{2}} P(L(\frac{1}{\sqrt{1+\delta^{2}}}, \frac{\delta}{\sqrt{1+\delta^{2}}}))  \textit{.}
\end{equation}
We know that $g_{L}$ is analytic because $\gamma$ is analytic. So, let's suppose that $f_{L}$ is equal to a polynomial. Let's observe that $ P(L(\frac{1}{\sqrt{1+\delta^{2}}}, \frac{\delta}{\sqrt{1+\delta^{2}}}))$ is bounded so that $g_{L} $ can only be of degree at most one. \\
From this fact, one has that, for every $\lambda, \delta \in \mathbb{R}$, $(P \circ L)(\lambda,\delta)$ can be written as : 
\begin{equation}
\label{eq1000020}
 (P \circ L)(\lambda,\delta) = a_{0}\lambda + a_{1} \delta
\end{equation}
where $a_{0},a_{1} \in \mathbb{R}$. \\
By making the change of variable $u = L(\lambda,\delta)$, one gets finally that there exists $b_{0},b_{1} \in \mathbb{R}$ such that for all $\lambda, \delta \in \mathbb{R}$ : 
\begin{equation}
\label{eq1000022}
P(\lambda,\delta) = b_{0} \lambda + b_{1} \delta \textit{.}
\end{equation}
Yet, by using Equation ($\ref{eq1000021}$), one gets that : 
\begin{equation}
\label{eq1000023}
x \partial_{x}Y_{\gamma} - y \partial_{y} Y_{\gamma} = b_{0} x + b_{1} y \textit{.}
\end{equation}
Yet, $Y_{\gamma}$ is positively homogeneous and so for $x > 0$, one has $Y_{\gamma}(x,y)= x Y_{\gamma}(1,\frac{y}{x})$. Using differentiation, one gets that : 
\begin{equation}
\label{eq1000024}
\partial_{y}Y_{\gamma}(x,y) = \partial_{y}Y_{\gamma} (1, \frac{y}{x}) \textit{.}
\end{equation}
Using Equation $(\ref{eq1000023})$ and Equation $(\ref{eq1000024})$, by diving by $x$ and making $x$ goes to infinity one gets that : 
\begin{equation}
\label{eq1000025}
\partial_{x}Y_{\gamma} (x,y) \rightarrow b_{0}
\end{equation} 
when $x \rightarrow \infty$. \\
So, one gets that (because $Y_{\gamma}$ is analytic) : 
\begin{equation}
\label{eq1000026}
\partial_{x}Y_{\gamma} (x,y) = K(y) 
\end{equation}
where $K(y)$ is a regular function. \\
By interchanging the role of $x$ and $y$, one gets that : 
\begin{equation}
\label{eq1000027}
Y_{\gamma}(x,y) = A_{0}xy + A_{1}
\end{equation}
where $A_{0}$ and $A_{1}$ belong to $\mathbb{R}$. Yet by setting $x=1$ and $y = \delta$ and by using Lemma 5.2 from $\cite{dolgopyat2014deviations}$, we obtain that it is impossible.
\end{proof}
In the non symmetric case, we will need the following lemma. 
\begin{lemma}
The following alternative holds. Let $k \geqslant 2$, $k \in 2 \mathbb{N}$. Either \\
(i) There exists $L \in GL_{2}(\mathbb{R})$ and $\delta,\delta' \in \mathbb{R}$ such that 
\begin{equation}
\label{eq1000028}
\frac{g_{L}^{(k)}(\delta)}{g_{L}^{(k)}(\delta')} \neq \frac{\tilde{g}_{L}^{(k)}(\delta)}{\tilde{g}_{L}^{(k)}(\delta')}
\end{equation}
or  \\
(ii) $\Omega_{\gamma}$ has a center of symmetry.
\end{lemma}
\begin{proof}
Let's suppose that (i) does not hold. Let $L = Id$. We have that $g_{Id}^{(k)}(\cdot) = c \tilde{g}_{Id}^{(k)}(\cdot)$ for some constant $c$. In other words
\begin{equation}
\label{eq1000029}
(\frac{\partial}{\partial \delta})^{k} P (1,\delta) = c (\frac{\partial}{\partial \delta})^{k} P (-1,-\delta) \textit{.}
\end{equation}
Since for $x > 0$ we have $P(x,y) = x P(1, \frac{y}{x})$, it follows that 
\begin{equation}
\label{eq1000030}
\partial_{y}^{k} P (x,y) = c \partial_{y}^{k} P (-x,-y) \textit{.}
\end{equation}
Since $\gamma$ is analytic, this equality in fact holds identically. In particular : 
\begin{equation}
\label{eq1000031}
\partial_{y}^{k} P (-x,-y) = c \partial_{y}^{k} P (x,y) \textit{.}
\end{equation}
As a consequence, and because of Lemma $\ref{lemme102}$, one must have $c= \pm 1$. \\
Furthermore, from Equation $(\ref{eq1000030})$, one has necessarily that : 
\begin{equation}
\label{eq1000032}
P(x,y) - c P(-x,-y) = \sum_{l=0}^{k-1} a_{l}(x) y^{l} \textit{.}
\end{equation}
Let's set now $H(x,y) = Y_{\gamma}(x,y) - c Y_{\gamma}(-x,-y)$. Equation $(\ref{eq1000032})$ and Equation ($\ref{eq1000032}$) give us that : 
\begin{equation}
\label{eq1000033}
x \partial_{x}H(x,y) - y \partial_{y}H(x,y) = \sum_{l=0}^{k-1} a_{l}(x) y^{l} \textit{.}
\end{equation}
Furthermore, one has, by positive homogeneity of $H(x,y)$, that $\partial_{x}H(x,y) = \partial_{x}H(x/y,1)$. \\
Equation $(\ref{eq1000033})$ gives us then that the growth of $\partial_{y}H(x,y)$, at $x$ fixed, is of order $y^{d-1}$. As $H(x,y)$ is analytic, we must have : 
\begin{equation}
\label{eq1000034}
\partial_{y}H(x,y) = \sum_{l=0}^{k-1} \tilde{b}_{l}(x) y^{l}
\end{equation}
and so 
\begin{equation}
\label{eq1000035}
H(x,y) = \sum_{l=0}^{k} b_{l}(x) y^{l}
\end{equation}
where the $b_{l}(x)$ are analytic functions. \\
Yet, $H(x,y)= y H(\frac{x}{y},1)$. So, one must have : $H(x,y) = b_{0}(x) + b_{1}(x)y \textit{.}$ 
So, we see that we have reached the same point of the demonstration of Lemma 5.3 of $\cite{dolgopyat2014deviations}$. So, by reasoning the same way, one gets finally that : 
\begin{equation}
\label{eq1000036}
Y_{\gamma}(X) - Y_{\gamma}(-X) = <X,v>
\end{equation}
with $v$ being a vector of $\mathbb{R}^{2}$. By shifting the origin to $x_{0}$, $Y_{\gamma}(X)$ is replaced by $Y_{\gamma}(X) + <X,x_{0}>$ and $Y_{\gamma}(-X)$ by $Y_{\gamma}(-X) - <X,x_{0}>$. So, after shifting the origin to $\frac{v}{2}$, we get $Y_{\gamma}(X) = Y_{\gamma}(-X)$ so that $\Omega_{\gamma}$ is symmetric.
\end{proof}
Let's assume WLOG that Equation ($\ref{eq1000028}$) holds for $L = Id$.
\begin{lemma}
\label{lemme101}
For every $m \in \mathbb{N}-\{ 0 \}$, for every family $(p_{1},\cdots,p_{m}) \in \mathbb{Z}^{m}$, for every $k_{1},\cdots,k_{m} \in \Pi$ all distinct, one has the following implication : 
if $ \sum_{i=1}^{m}  p_{i} f_{k_{i}} $ is a polynomial function then all the $p_{i}$ must be equal to $0$. \\
In the non symmetric case, for every $m \in \mathbb{N}-\{ 0 \}$, for every family $(p_{-m},\cdots,p_{-1},p_{1},\cdots,p_{m}) \in \mathbb{Z}^{2m}$, for every $k_{1},\cdots,k_{m} \in \Pi$ all distinct if 
$\sum_{i=1}^{m}  p_{i} f_{k_{i}} + \sum_{i=1}^{m}  p_{-i} f_{-k_{i}} $ is a polynomial function then all the $p_{i}$ must be equal to $0$.
\end{lemma}
\begin{proof}
Let's suppose that $ \sum_{i=1}^{m}  p_{i} f_{k_{i}} $ is a polynomial. Let's give us $j \in [1,m]$. We are going to show that $p_{j} = 0$. \\
Let $\beta \in \mathbb{R}^{2}$ such that $<k_{j},\beta> \neq 0$. Let $\alpha \in \mathbb{R}^{2}$ and let $X= \alpha$ and $Y= \delta \alpha + \theta \beta$. Then, one has : $$f_{k}(A,B) = |<k,\alpha>|g_{Id}(\delta + \theta \frac{<k,\beta>}{|<k,\alpha>|}) $$ if $<k,\alpha> > 0$
or
$$f_{k}(A,B) = |<k,\alpha>|\tilde{g}_{Id}(\delta + \theta \frac{<k,\beta>}{|<k,\alpha>|}) $$
if $<k,\alpha> < 0$. \\
As $ \sum_{i=1}^{m}  p_{i} f_{k_{i}} $ is a polynomial then $ \sum_{i=1}^{m}  p_{i} f_{k}(A,B) $ is a polynomial in $\beta$ whose degree is bounded by a number that does not depend on $\alpha$ (nor in $\beta$). \\
So there exists $ K \geqslant 2$ such that the $K$-th derivative of $ \sum_{i=1}^{m} f_{k}(A,B) $ relatively to $\theta$ is equal to $0$. It means that the terms in front of $\theta^{K}$ is equal to $0$. Hence the following equation : 
\begin{equation}
\label{eq1007}
\sum_{i=1}^{m}  h_{i} \frac{<k_{i},\beta>^{K}}{|<k_{i},\alpha>|^{K-1}} = 0 
\end{equation}
where $h_{i} = p_{i} g_{Id}^{(K)}(\delta)$ if $<k_{i},\alpha> > 0$ and $h_{i} = p_{i} g_{Id}^{(K)}(\delta)$ if $<k_{i},\alpha> < 0$.
Now, since the $k_{i}$ belong to $\Pi$ and are all distinct, it is possible to choose $\alpha$ so that $<k_{j},\alpha> > 0$ is arbitrary small while $<k_{i},\alpha> $ remain bounded away from zero for every $i \neq j$. Thus, we must have $h_{j} = 0$. Yet, $g$ is not a polynomial so there exists $\delta \in \mathbb{R}$ such that $g_{Id}^{(K)}(\delta) \neq 0$ and it gives us that $p_{j} = 0$. 
In the non symmetric case, Equation $(\ref{eq1007})$ becomes 
\begin{equation}
\label{eq1000037} 
\sum_{i=1}^{m}  h_{i} \frac{<k_{i},\beta>^{K}}{|<k_{i},\alpha>|^{K-1}} = 0 \textit{.}
\end{equation}
where $h_{i} = p_{i} g_{Id}^{(K)}(\delta) + p_{-i} \tilde{g}_{Id}^{(K)}(\delta)$ if $<k_{i},\alpha> > 0$ and  $h_{i} = p_{i} \tilde{g}_{Id}^{(K)}(\delta) + p_{-i} g_{Id}^{(K)}(\delta)$ if $<k_{i},\alpha> < 0$. \\
Let's consider for example the case where the first alternative holds. As before, we must have $p_{j} g_{Id}^{(K)}(\delta) + p_{-j} \tilde{g}_{Id}^{(K)}(\delta) = h_{j} = 0$ for any choice of $\delta$. Since we assume that Equation $(\ref{eq1000028})$, this implies that $ p_{j} = p_{-j} = 0$.
\end{proof}
Lemma $\ref{lemme101}$ enables us to prove the following lemma.
\begin{lemma}
\label{lemme100}
For a typical $L \in \mathscr{S}_{2}$, for every $m \in \mathbb{N}-\{ 0 \}$, for every family $(p_{1},\cdots,p_{m}) \in \mathbb{Z}^{m}$, for every $k_{1},\cdots,k_{m} \in \Pi$ all distinct if 
\begin{equation}
\label{eq505}
\sum_{i=1}^{m}  p_{i} W_{k_{i}}(L) = 0
\end{equation}
then $$\sum_{i=1}^{m}  p_{i} f_{k_{i}} = 0 \textit{.} $$ 
In the non symmetric case, for a typical $L \in \mathscr{S}_{2}$, for every $m \in \mathbb{N}-\{ 0 \}$, for every family $(p_{-m},\cdots,p_{-1},p_{1},\cdots,p_{m}) \in \mathbb{Z}^{2m}$, for every $k_{1},\cdots,k_{m} \in \Pi$ all distinct if 
\begin{equation}
\label{eq0505}
\sum_{i=1}^{m}  p_{i} W_{k_{i}}(L) + \sum_{i=1}^{m} p_{-i} W_{-k_{i}}(L) = 0 
\end{equation}
then
$$ \sum_{i=1}^{m}  p_{i} f_{k_{i}} + \sum_{i=1}^{m} p_{-i}f_{-k_{i}} = 0 \textit{.}$$
\end{lemma}
\begin{proof}
First, we note that, for every $k \in -\Pi \cupdot \Pi$, $W_{k}(L) = f_{k}((e_{1}(L))_{1},(e_{2}(L))_{1},(e_{1}(L))_{2},(e_{2}(L))_{2}) \textit{.}$
Thanks to this remark, we see that Lemma $\ref{lemme100}$ is a direct consequence of the facts that $\mathscr{S}_{2} = SL_{2}(\mathbb{R}) / SL_{2}(\mathbb{Z})$, that $\det$ is a polynomial function whereas $ \sum_{i=1}^{m}  p_{i} f_{k_{i}} $ and $ \sum_{i=1}^{m}  p_{i} f_{k_{i}} + \sum_{i=1}^{m} p_{-i}f_{-k_{i}}$ are polynomial functions, for a typical $L$, only in the trivial case according to Lemma $\ref{lemme101}$.
\end{proof}
We can now give the proof of Proposition $\ref{prop13}$.
\begin{proof}[Proof of Proposition $\ref{prop13}$]
It is a direct consequence of Lemma $\ref{lemme101}$ and Lemma $\ref{lemme100}$.
\end{proof}
\subsection{Proof of Proposition $\ref{prop12}$}
As we have at our disposal Proposition $\ref{prop13}$, the proof of Proposition $\ref{prop12}$ is the same as the proof of Proposition 5 in $\cite{trevisan2021limit}$. We are going to present it again for completeness. \\
Before starting the proof of Proposition $\ref{prop12}$, we only need a simple lemma : 
\begin{lemma}
\label{lemme19}
For every $m \in \mathbb{N}-\{ 0 \}$, for every family $(p_{1},\cdots,p_{m}) \in \mathbb{Z}^{m}-\{ 0 \}$, for every $k_{1},\cdots,k_{m} \in \Pi$, all distinct, for every $0 < \epsilon < 1$, there exists $a > 0$ and $K_{\epsilon}$ a measurable set of $\mathscr{S}_{2}$ such that $\tilde{\mu}_{2} (\mathscr{S}_{2}-K_{\epsilon}) \leqslant \epsilon$ and such that for all $L \in  K_{\epsilon}$,
$$ \lvert \sum_{i=1}^{m}  p_{i} W_{k_{i}}(L) \rvert \geqslant a \textit{ .} $$
In the non symmetric case, for every $m \in \mathbb{N}-\{ 0 \}$, for every family $(p_{-m},\cdots,p_{-1},p_{1},\cdots,p_{m}) \in \mathbb{Z}^{2m}-\{ 0 \}$, for every $k_{1},\cdots,k_{m} \in \Pi$, all distinct, for every $0 < \epsilon < 1$, there exists $a > 0$ and $K_{\epsilon}$ a measurable set of $\mathscr{S}_{2}$ such that $\tilde{\mu}_{2} (\mathscr{S}_{2}-K_{\epsilon}) \leqslant \epsilon$ and such that for all $L \in  K_{\epsilon}$,
$$ \lvert \sum_{i=1}^{m}  p_{i} W_{k_{i}}(L) + \sum_{i=1}^{m}  p_{-i} W_{-k_{i}}(L)  \rvert   \geqslant a \textit{ .} $$
\end{lemma}
\begin{proof}
It is a direct consequence of Proposition $\ref{prop13}$.
\end{proof}
Now, by using the foliation given by $\delta(\lambda)$ and previous results, we can now prove Proposition $\ref{prop12}$.
\begin{proof}[Proof of Proposition $\ref{prop12}$]
The proof in all its generality can be made as in the case where $\tilde{\mu}_{2} = \mu_{2}$. So, we will suppose for simplicity that $\tilde{\mu}_{2} = \mu_{2}$. \\
We consider the case where $\Omega_{\gamma}$ is symmetric. The non symmetric case is similar. \\
Let $l \geqslant 1$. Let $\psi \in C_{c}^{\infty}(\mathscr{S}_{2})$. Let $(p_{1},\cdots,p_{l}) \in \mathbb{Z}^{l}-\{ 0 \}$. \\
For all $\epsilon > 0$, we call $\mathcal{F}_{\epsilon}$ the tribe on $\mathscr{S}_{2}$ generated by the $J_{\epsilon}(L)$. Let $1 > \epsilon_{1} \geqslant \epsilon_{2} > 0$. \\
According to Lemma $\ref{lemme11}$ and Lemma $\ref{lemme19}$, there exists a measurable part $K_{\epsilon_{1}}$ such that $\mu_{2}(K_{\epsilon_{1}}) \geqslant 1- \epsilon_{1}$, a real $M > 0$ and a real $a > 0$, such that
\begin{itemize}
\item for every $L \in K_{\epsilon_{1}}$, for every $\lambda \in [\frac{1}{1+\epsilon_{1}}, 1+ \epsilon_{1}]$ : 
\begin{equation}
\label{eq87}
|\psi(\delta(h)L) -\psi(L)|  \leqslant M |h| 
\end{equation}
where $h= \lambda -1$
\item for every $L \in K_{\epsilon_{1}}$, for every $\lambda \in [\frac{1}{1+\epsilon_{1}}, 1+ \epsilon_{1}]$, Equation ($\ref{eq61}$) is verified.
\item for every $L \in K_{\epsilon_{1}}$, 
\begin{equation}
\label{eq300}
\lvert \sum_{i=1}^{l}  p_{i} W_{k_{i}}(L) \rvert \geqslant a \textit{ .} 
\end{equation}
\end{itemize}
Furthermore, we are going to suppose, even if it means making $\epsilon_{2}$ goes to $0$, that for every $L \in K_{\epsilon_{1}}$, 
\begin{equation}
\label{eq506}
J_{\epsilon_{2}}(L) \subset K_{\epsilon_{1}} \textit{.}
\end{equation}
$\textbf{Claim. }$With these notations, we have, for all $L \in K_{\epsilon_{1}}$, that : 
\begin{equation}
\label{eq93}
 \mathbb{E} \left( \psi e(\sum_{i=1}^{l} p_{i} \theta_{k_{i}}) | \mathcal{F}_{\epsilon_{2}} \right) (L) = O(\epsilon_{2} ) + O(\frac{1}{a t \epsilon_{2}}) + \frac{1}{a} O(\epsilon_{2}) \textit{.}
\end{equation}
To this end,  let's set, for all $ \epsilon > 0$, $$ \delta(\epsilon) = 1+\epsilon - \frac{1}{1+\epsilon} \textit{.}$$ \\
Then one has for every $L \in K_{\epsilon_{1}}$ according to Lemma $\ref{lemme11}$ :
\begin{align}
\label{eq89}
& \mathbb{E}(\psi e(\sum_{j=1}^{l} p_{j} \theta_{k_{j}}) | \mathcal{F}_{\epsilon_{2}} )(L) \nonumber \\
& = \frac{1}{\delta(\epsilon_{2})} \int_{\frac{1}{1+\epsilon_{2}}-1}^{\epsilon_{2}} \left( \psi e(\sum_{j=1}^{l} p_{j} \theta_{k_{j}}) \right) (\delta(h) L) d h \nonumber \\
& = \psi(L) \frac{1}{\delta(\epsilon_{2})} \int_{\frac{1}{1+\epsilon_{2}}-1}^{\epsilon_{2}} \left(e(\sum_{j=1}^{l} p_{j} \theta_{k_{j}}) \right) (\delta(h) L) d h + O(\epsilon_{2}) \nonumber \\
& = \left( \psi e(\sum_{j=1}^{l} p_{j} \theta_{k_{j}}) \right)(L)  \frac{1}{\delta(\epsilon_{2})} \int_{\frac{1}{1+\epsilon_{2}}-1}^{\epsilon_{2}} e^{i t D_{1}(L) h + i t D_{2}(L,h) } d h + O(\epsilon_{2})
\end{align}
where
\begin{equation}
\label{eq90}
D_{1}(L) = \sum_{j=1}^{l} p_{j} W_{k_{j}}(L) \textit{, }
\end{equation}
\begin{equation}
\label{eq91}
D_{2}(L,h) = \sum_{j=1}^{l} p_{j} \theta_{k_{j}}(\delta(\lambda) L) -\sum_{j=1}^{l} p_{j} \theta_{k_{j}}(L)- D_{1}(L)h  \textit{ such that }
\end{equation}
\begin{equation}
\label{eq92}
D_{2}(L,h) = O(h^{2}) \textit{ and }
\end{equation}
$D_{2}(L, \cdot)$ is smooth around $0$. \\
Thus, by integrating by part and by using Equation ($\ref{eq300}$), one gets that for all $L \in K_{\epsilon_{1}}$ : 
\begin{align}
\label{eq507}
& \frac{1}{\delta(\epsilon_{2})} \int_{\frac{1}{1+\epsilon_{2}}-1}^{\epsilon_{2}} e^{i t D_{1}(L) h + i t D_{2}(L,h) } d h \nonumber \\
& = \frac{1}{\delta(\epsilon_{2})} \left( \left[  \frac{e^{i t D_{1}(L) h + i t D_{2}(L,h)}}{i t D_{1}(L) } \right]_{\frac{1}{1+\epsilon_{2}}-1}^{\epsilon_{2}} + \frac{1}{D_{1}(L)} \int_{\frac{1}{1+\epsilon_{2}}-1}^{\epsilon_{2}} \left( D_{2}(L,\cdot) \right) '(h) e^{i t D_{1}(L) h + i t D_{2}(L,h)} dh \right) \nonumber \\
& = O(\frac{1}{a t \epsilon_{2}}) + \frac{1}{a} O(\epsilon_{2}) \textit{.}
\end{align}
Finally, Equation ($\ref{eq89}$) and Equation ($\ref{eq507}$) give the wanted claim. \\
\\
Thanks to Equation ($\ref{eq93}$), the fact that $\mu_{2}(K_{\epsilon_{1}}) \geqslant 1- \epsilon_{1}$ and because of Equation ($\ref{eq506}$), we have that : 
\begin{align}
\label{eq94}
& |\mathbb{E}\left(\psi e(\sum_{j=1}^{l} p_{j} \theta_{k_{j}}) \right) | \nonumber \\
& \leqslant  |\mathbb{E} \left( \psi e(\sum_{j=1}^{l} p_{j} \theta_{k_{j}}) \mathbf{1}_{K_{\epsilon_{1}}^{c}} \right)|  +  |\mathbb{E}\left( \psi e(\sum_{j=1}^{l} p_{j} \theta_{k_{j}}) \mathbf{1}_{K_{\epsilon_{1}}} \right) | \nonumber \\
& \leqslant \lVert \psi \rVert_{\infty} \epsilon_{1} + O(\epsilon_{2} ) + O(\frac{1}{a t \epsilon_{2}}) + \frac{1}{a} O(\epsilon_{2}) \textit{.}
\end{align}
 By first choosing $\epsilon_{1} > 0$ small enough (note that $a$ depends on $\epsilon_{1}$), then choosing $\epsilon_{2} > 0$ and finally choosing $t$ large enough, we obtain the wanted result. 
\end{proof}
We are now brought back to the study of the convergence, in the symmetric case, of the sums of the type $$\tilde{S}_{A}(\omega,L) = \sum_{k \in \Pi_{A}(L)} \frac{\rho_{\gamma}(k_{1} e_{1}(L) + k_{2} e_{2}(L)) \tilde{Z}_{k}(\omega)}{\lVert k_{1} e_{1}(L) + k_{2} e_{2}(L) \rVert^{\frac{3}{2}}}   $$
and of the sums of the type, in the non symmetric case,
$$\tilde{S}_{A}(\omega_{1},\omega_{2},L) = \sum_{k \in \Pi_{A}(L)} \frac{\rho_{\gamma}(k_{1} e_{1}(L) + k_{2} e_{2}(L)) \tilde{Z}_{k}(\omega_{1} ) + \rho_{\gamma}(-(k_{1} e_{1}(L) + k_{2} e_{2}(L))) \tilde{Z}_{k}(\omega_{2} )}{\lVert k_{1} e_{1}(L) + k_{2} e_{2}(L) \rVert^{\frac{3}{2}}}   $$
where $Z_{k}$ are non-zero real independent identically distributed random variables from $\Omega \ni \omega$ are symmetrical and have a compact support and where $(\omega_{1}, \omega_{2}) \in \Omega \times \Omega$. In the next section, we are going to study the sums of this type and prove Proposition $\ref{prop15}$, which will conclude the proof of Theorem $\ref{thm2}$. 
\section{Asymptotic study of $\tilde{S}_{A}(\omega,L)$ and of $\tilde{S}_{A}(\omega_{1},\omega_{2},L)$} 
The study is in fact essentially done in Section 5 of $\cite{trevisan2021limit}$. We only need to notice the following fact to apply the same method : 
\begin{lemma}
\label{lemme103}
There exists $m,M > 0$ such that for every $\xi \in \mathbb{R}^{2}-\{0\}$, $$M \geqslant \rho_{\gamma}(\xi) \geqslant m \textit{.} $$
\end{lemma}
\begin{proof}
As $\gamma$ is strictly convex, from the definition of $\rho_{\gamma}$, $\rho_{\gamma} > 0$. \\
Furthermore, for every $\xi \in \mathbb{R}^{2}-\{0 \}$, $$\rho_{\gamma}(\xi) = \rho_{\gamma}(\frac{\xi}{\lVert \xi \rVert}) \textit{.}$$
$S^{1}$ being compact, we obtain immediately the wanted result.
\end{proof}
\begin{proof}[Proof of Proposition $\ref{prop15}$]
We use the same method of the last section of $\cite{trevisan2021limit}$ and use Lemma $\ref{lemme103}$.
\end{proof}
\section{Elements of proof of Theorem $\ref{thm4}$}
To prove Theorem $\ref{thm4}$ we follow the same approach that was followed to prove Theorem $\ref{thm3}$. \\
Namely, as a first step, we prove the equivalent of Proposition $\ref{prop8}$. Namely, let's set 
\begin{equation}
\label{eq1000038}
\Delta_{A,prime}(L,\alpha,t) = |\frac{\mathcal{R}(t \Omega_{\gamma} + \alpha,L)}{\sqrt{t}} - S_{A,prime}(L,\alpha,t)| 
\end{equation}
with
\begin{align}
\label{eq1000039}
& S_{A,prime}(L,\alpha,t) = \frac{1}{\pi} \sum_{\substack{ l \in L^{\perp} \textit{ prime}  \\ 0 < \lVert l \rVert \leqslant A }} \frac{ 1}{\lVert l \rVert^{\frac{3}{2}}}  \\ 
& \sum_{m \in \mathbb{N}-\{ 0 \}} \frac{ \rho_{\gamma}(l) \cos(2 \pi t m Y_{\gamma}(l) + 2 \pi m <\alpha,l> - \frac{3 \pi}{4}) + \rho_{\gamma}(-l) \cos(2 \pi t m Y_{\gamma}(-l) - 2 \pi m <\alpha,l> - \frac{3 \pi}{4})  }{m^{\frac{3}{2}}}  \nonumber \textit{.}
\end{align}
Then, one can prove, as Proposition $\ref{prop8}$ was proven : 
\begin{prop}
\label{prop20}
For every $\beta > 0$, for every $A > 0$ large enough, for every $t$ large enough, one has that :
$$\mathbb{P}(\Delta_{A,prime}(L,\alpha,t)  \geqslant \beta) \leqslant \beta \textit{.}$$
\end{prop}
Second, as Proposition $\ref{prop11}$ was already proven, the limit distribution of $S_{A,prime}(L,\alpha,t)$ when $t \rightarrow \infty$ and with $\alpha \in \mathbb{R}^{2}$ fixed and with $L$ being distributed according to $\tilde{\mu}_{2}$ is, in the symmetric case, is the distribution of the almost-sure limit of 
\begin{equation}
\label{eq1000040}
S_{\gamma,A}(\theta,L) = \frac{2}{\pi} \sum_{\substack{ (k_{1},k_{2}) \in \Pi  \\ \lVert e(k_{1},k_{2},L) \rVert \leqslant A }} \frac{\rho_{\gamma}(e(k_{1},k_{2},L)) \phi_{\alpha}(\theta_{(k_{1},k_{2})},e(k_{1},k_{2},L))}{\lVert e(k_{1},k_{2},L) \rVert^{\frac{3}{2}}}
\end{equation}
when $A \rightarrow \infty$ and where $\theta = (\theta_{(k_{1},k_{2})}) \in \mathbb{T}^{\infty}$ being distributed according to $\lambda_{\infty}$ and $L$ being distributed according to $\tilde{\mu}_{2}$ and where $e(k_{1},k_{2},L) = k_{1} e_{1}(L) + k_{2} e_{2}(L)$. \\
In the non symmetric case, the limit distribution of $S_{A,prime}(L,\alpha,t)$ when $t \rightarrow \infty$ and with $\alpha \in \mathbb{R}^{2}$ fixed and with $L$ being distributed according to $\tilde{\mu}_{2}$ is the distribution of the almost-sure limit of 
\begin{equation}
\label{eq1000041}
S_{\gamma,A}(\theta,L) = \frac{1}{\pi} \sum_{\substack{ (k_{1},k_{2}) \in \Pi  \\ \lVert e(k_{1},k_{2},L) \rVert \leqslant A }} \frac{ \phi_{\alpha,\gamma,2}(\theta_{(k_{1},k_{2})},e(k_{1},k_{2},L))}{\lVert e(k_{1},k_{2},L) \rVert^{\frac{3}{2}}}
\end{equation}
when $A \rightarrow \infty$ and where $\theta = (\theta_{(k_{1},k_{2})}) = (\theta_{1,(k_{1},k_{2})},\theta_{2,(k_{1},k_{2})}) \in \mathbb{T}^{\infty,2}$ being distributed according to $\lambda_{\infty,2}$ and $L$ being distributed according to $\tilde{\mu}_{2}$. \\
Finally, we conclude the proof by proving the equivalent of Proposition $\ref{prop15}$. To do so, we follow the exact same approach that was used to prove Proposition $\ref{prop15}$.

\bibliographystyle{plain}
\bibliography{bibliographie}

\begin{thebibliography}{10}

\bibitem{gorodnikcounting}
Michael Björklund and Alexander Gorodnik.
\newblock Counting in generic lattices and higher rank actions, 2021.

\bibitem{bleher1992distribution}
Pavel Bleher et~al.
\newblock On the distribution of the number of lattice points inside a family
  of convex ovals.
\newblock {\em Duke Mathematical Journal}, 67(3):461--481, 1992.

\bibitem{bleher1993distribution}
Pavel~M Bleher, Zheming Cheng, Freeman~J Dyson, and Joel~L Lebowitz.
\newblock Distribution of the error term for the number of lattice points
  inside a shifted circle.
\newblock {\em Communications in mathematical physics}, 154(3):433--469, 1993.

\bibitem{bassam}
Dmitry Dolgopyat and Bassam Fayad.
\newblock Deviations of ergodic sums for toral translations ii. boxes.
\newblock {\em arXiv preprint arXiv:1211.4323}, 2012.

\bibitem{dolgopyat2014deviations}
Dmitry Dolgopyat and Bassam Fayad.
\newblock Deviations of ergodic sums for toral translations i. convex bodies.
\newblock {\em Geometric and Functional Analysis}, 24(1):85--115, 2014.

\bibitem{hardy1917average}
GH~Hardy.
\newblock The average order of the arithmetical functions p (x) and $\delta$
  (x).
\newblock {\em Proceedings of the London Mathematical Society}, 2(1):192--213,
  1917.

\bibitem{heath1992distribution}
DR~Heath-Brown.
\newblock The distribution and moments of the error term in the dirichlet
  divisor problem.
\newblock 1992.

\bibitem{heath2012lattice}
DR~Heath-Brown.
\newblock Lattice points in the sphere.
\newblock In {\em Number theory in progress}, pages 883--892. de Gruyter, 2012.

\bibitem{huxley2003exponential}
Martin~N Huxley.
\newblock Exponential sums and lattice points iii.
\newblock {\em Proceedings of the London Mathematical Society}, 87(3):591--609,
  2003.

\bibitem{iwaniec1988divisor}
Henryk Iwaniec and CJ~Mozzochi.
\newblock On the divisor and circle problems.
\newblock {\em Journal of Number theory}, 29(1):60--93, 1988.

\bibitem{kelmer2021second}
Dubi Kelmer and Shucheng Yu.
\newblock The second moment of the siegel transform in the space of symplectic
  lattices.
\newblock {\em International Mathematics Research Notices}, 2021(8):5825--5859,
  2021.

\bibitem{kesten60u}
Harry Kesten.
\newblock Uniform distribution mod 1.
\newblock {\em Annals of Mathematics}, pages 445--471, 1960.

\bibitem{kesten62u}
Harry Kesten.
\newblock Uniform distribution mod 1 (ii).
\newblock {\em Acta Arithmetica}, 4(7):355--380, 1962.

\bibitem{marklof1998n}
Jens Marklof and Zeev Rudnick.
\newblock The n-point correlations between values of a linear form.
\newblock Technical report, SCAN-9905042, 1998.

\bibitem{trevisan2021limit}
Julien Trevisan.
\newblock Limit laws in the lattice problem. ii. the case of ovals, 2021.

\bibitem{vinogradov2010limiting}
Ilya Vinogradov.
\newblock Limiting distribution of visits of sereval rotations to shrinking
  intervals.
\newblock {\em arXiv preprint arXiv:1005.1622}, 2010.

\end{thebibliography}
\end{document}